\documentclass[11pt,reqno]{amsart}
\usepackage{epsfig}
\numberwithin{equation}{section}

\newtheorem{theorem}{Theorem}[section]

\newtheorem{proposition}[theorem]{Proposition}

\newtheorem{remark}[theorem]{Remark}
\usepackage{epsfig}
\usepackage{psfrag}
\usepackage{amssymb}
\usepackage{amsthm}
\usepackage{amsmath}
\usepackage{graphicx}

 \newtheorem{thm}{Theorem}[section]
 \newtheorem{cor}[thm]{Corollary}

 \theoremstyle{definition}
 \newtheorem{defn}[thm]{Definition}
 \newtheorem{rem}[thm]{Remark}
 \numberwithin{equation}{section}
 \def\al{\aligned}
\def\eal{\endaligned}
\def\be{\begin{equation}}
	\def\ee{\end{equation}}
\def\lab{\label}
\def\td{\tilde}
\def\af{\alpha}

\def\al{\aligned}
\def\eal{\endaligned}
\def\be{\begin{equation}}
\def\ee{\end{equation}}
\def\lab{\label}

\def\al{\aligned}

\numberwithin{equation}{section}

\usepackage[backref=page]{hyperref}

\begin{document}

\title[]{Dimension reduction of axially symmetric Euler equations near  maximal points off the axis}
\author{Qi S. Zhang}
\address{Department of Mathematics, University of California, Riverside, CA 92521, USA}
\email{qizhang@math.ucr.edu}
\date{October 05 2023; MSC2020: 35Q31}

\begin{abstract}
Let $v$ be a solution of the axially symmetric Euler
equations (ASE) in a finite cylinder in $\mathbb{R}^3$. We show that suitable blow-up limits of possible velocity singularity and most self similar vorticity singularity near maximal points off the vertical axis are  two dimensional ancient solutions of the Euler equation in either $\mathbb{R}^2 \times (-\infty, 0]$ or $\mathbb{R}^2_+ \times (-\infty, 0]$.
This reduces the search of off-axis self-similar or other velocity blow-up solutions to a problem involving purely 2-dimensional Euler equations.
 Also, some asymptotic self-similar velocity blow-up and  expected asymptotic self-similar vorticity blow up scenario  at  the boundary appear to be ruled out. On the other hand, this method may provide a path to velocity blow up if one can construct certain stable ancient solutions to the 2-d Euler equation in the half plane.

\end{abstract}
\maketitle

\section{Introduction}

In this paper we take a step in the analytical study of the structure, in a space-time region with
high flow speed or high vorticity, of solutions to the three dimensional
incompressible Euler equations
\begin{equation}\label{EL}
\begin{cases}
\partial_t v + v\cdot \nabla v + \nabla p = 0,\\[-4mm]\\
    \nabla \cdot v = 0,
\end{cases}
\quad t \in [0, T), \quad x \in D \subset \mathbb{R}^3
\end{equation}
with  axially symmetric data. Namely,
in cylindrical coordinates, the solution $v=v(x, t)$ is of the
form
\begin{equation}\label{axi-solu}
v(x, t) = v^r(r, x^{(3)}, t)e_r +  v^\theta(r, x^{(3)}, t)e_\theta +  v^{(3)}(r,
x^{(3)}, t)e_3.
\end{equation}
Here $x = (x^{(1)}, x^{(2)}, x^{(3)})$, $r = \sqrt{(x^{(1)})^2 + (x^{(2)})^2}$ and
\begin{equation}
e_r =(\frac{x^{(1)}}{r},  \frac{x^{(2)}}{r},  0),
\quad e_\theta = (- \frac{x^{(2)}}{r}, \frac{x^{(1)}}{r}, 0),
\quad e_3 = (0, 0, 1)
\end{equation}
are the three orthogonal unit vectors along the radial,
angular, and axial directions respectively. It is well known, the
radial
, swirl and axial components $v^r$, $v^\theta$ and $v^{(3)}$
obey the equations
\begin{equation}\label{axi-EL}
\begin{cases}
\partial_tv^r + b\cdot\nabla v^r - \frac{(v^\theta)^2}{r}
  + \partial_rp = 0,\\[-4mm]\\
\partial_tv^\theta + b\cdot\nabla v^\theta + \frac{v^rv^\theta}{r}
  = 0,\\[-4mm]\\
\partial_t v^{(3)} + b\cdot\nabla v^{(3)} + \partial_{x^{(3)}} p =0,\\[-4mm]\\
b = v^re_r + v^{(3)} e_3,\quad \nabla\cdot b = \partial_r v^r +
\frac{v^r }{r} + \partial_{x^{(3)}} v^{(3)} = 0.
\end{cases}
\end{equation}

The axially symmetric case, studied a long time ago in \cite{UY} by  Ukhovskii and Yudovich e.g.,  appears much more special than the full
Euler equations.  However recent impressive research by several authors indicates that a regular solution with a certain initial value can form singularity in finite time even in this special case, especially in the presence of a boundary. In \cite{El} Elgindi proved a certain $C^{1, \alpha}$ initial data in $\mathbb{R}^3$ can generate finite time vorticity blow-up. See also \cite{EGM} by Elgindi, Ghoul and Masmoudi. In \cite{EJ} Elgindi and Jeong established finite time blow-up of some $C^{1, \alpha}$  solutions in the exterior of a cone. In the papers \cite{LH1, LH2} Luo and Hou provided strong numerical evidence that finite time blow-up of smooth solutions at boundary can occur. In two recent  pre-prints \cite{CH1, CH2}, Chen and Hou presented a computer assisted result indicating finite time vorticity blow-up of smooth solutions at the boundary of  a finite cylinder. See also \cite{WLG-SB} by Wang,  Lai,  Gomez-Serrano,  Buckmaster where a program of finding finite time self-similar blow-ups is initiated using neural networks. The common strategy is to identify certain asymptotic self-similar blow-up solutions and prove that they are stable under certain truncation and perturbation. Very often this is done by studying the associated  Boussinesq's equations. This seems to be one reason that some studies in the area are highly sophisticated and elongated.

In this paper, we prove that $C^{\gamma}_{x, t}$ blow-up limits (c.f. Definition \ref{defconv} )of all potential velocity singularity and most vorticity singularity of the axially symmetric Euler equations near maximal or anchor points off the rotational axis are ancient solutions to the  two dimensional Euler equations in $\mathbb{R}^2$ or $\mathbb{R}^2_+$, the half plane.
This  reduces the search of off-axis velocity blow-up and most asymptotic self similar vorticity blow up solutions to a problem involving purely
2-dimensional Euler equations.  Compared with the usual method of converting to the
2-dimensional Boussinesq's equations, which fits the study of the remaining vorticity blow up, this method seems to significantly reduce the complexity since the latter equations contain one extra unknown function, the temperature in addition to the velocity. In addition, certain asymptotic self-similar vorticity blow-ups at the boundary are apparently ruled out, c.f. Proposition \ref{prLH} (ii), \ref{prLH00} and \ref{prLH0} and the remarks below. The dimension reduction technique was previously used in \cite{LZ} to study axially symmetric Navier-Stokes equations.

Let us introduce some concepts and notations before stating the main result.

\begin{defn}
\lab{defbaopuops}

 Let $v$ be a smooth solution to \eqref{EL} in the domain $D \times [0, T)$ where $D \subset \mathbb{R}^3$ is an open set and $T>0$. Suppose the solution $v$ forms a singularity at the point $(x_*, T) \in D \times [0, T]$.

(a).
  A  sequence of points $\{ (x_k, t_k) \} \subset D \times [0, T)$  such that $x_k \to x_*$ and $t_k \to T$ as $k \to \infty$ are called blow-up points if $\lim_{k \to \infty} |v(x_k, t_k)|= \infty$. The points $x_k$ are called {\it off the  axis} if $|x'_k|>d>0$ for a fixed positive number $d>0$. Here $x' = (x^{(1)}, x^{(2)}, 0)$ if $x=(x^{(1)}, x^{(2)}, x^{(3)})$.

(b). Let $\{ Q_k \}$ be a sequence of positive numbers such that
\[
c^{-1} |v(x_k, t_k)| \le Q_k \le c |v(x_k, t_k)|, \quad k=1,2,3 ...,
\] for a fixed number $c \ge 1$.  The sequence
\[
\tilde v_k= \tilde v_k(\tilde x, \tilde t) = Q_k^{-1} v(Q_k^{-(1-\alpha)/\alpha} \tilde x + x_k,
Q_k^{-1/\alpha} \tilde t + t_k), \quad \tilde t \le 0, \quad \alpha \in (0, 1),
\]
is called a sequence of blow-up solutions with center $(x_k, t_k)$. Here and always $(\td x, \td t)$ are such that the unscaled variables $(x, t)$ are in the original domain.

(c).  The blow-up points $(x_k, t_k)$ are called anchor points if
\[
|\tilde v_k(\tilde x, \tilde t)| \le \lambda(|\td{x}|^2 +|\td t|)  |\tilde v_k(0, 0)|,
\quad \forall \quad (\td x, \td t) \quad \text{in domain such that } \quad  |\td{x}|^2 +|\td t| \le O(k),
\]for a continuous, positive, increasing function $\lambda: [0, \infty) \to [1, \infty)$.

(d). The blow-up points are called  {\it near maximal points} if
\[
Q_k=|v(x_k, t_k)| \ge c \sup_{s \in [0, t_k], \ y \in
D}  |v(y, s)|
\]for a fixed positive number $c>0$.

\end{defn}

It is clear that any near maximal  points are anchor points. The above definition can also be extended to the case of $\alpha<0$ and $Q_k \to 0$ as $k \to \infty$. In this case the velocity is scaled up together with the space time variables.

\begin{defn}
 \lab{defconv} Let $\{ v_k, p_k \}$ be a sequence of solutions to the Euler equation. We say it converges to another solution $\{v, p \}$ in local $C^\gamma_{x, t}$ sense if the sequences
$\{ v_k \}$, $\{ \nabla v_k \}$, $\{ \partial_t v_k \}$ and $\{ \nabla p_k \}$, converge in every compact set in the domain, to functions, $\{ v \}$, $\{ \nabla v \}$, $\{ \partial_t v \}$ and $\{ \nabla p \}$ respectively in $C^\gamma_{x, t}$ norm, i.e. the H\"older norm with exponent $\gamma$ and the parabolic distance $|x-y|+\sqrt{|t-s|}$.
Also the space $C^{2, 1, \gamma}_{x, t}$ is the usual one such that its elements $v$ together with  $\nabla v,  \nabla^2 v, \partial_t v$ are in $C^\gamma$ with respect to the parabolic distance. The space $C^{2, 1, \gamma}_{x, t, loc}$ is consisted of elements which are in $C^{2, 1, \gamma}_{x, t}$ in every compact set in the domain.
\end{defn}

\begin{defn}
\lab{defASSS}
A solution $v$ to \eqref{EL} is called  {\it asymptotic self-similar} (ASSS) with profile $V$ and error $W$ if
\[
v=v(x, t) = \frac{1}{(T_0-t)^\alpha} V \left( \frac{x}{(T_0-t)^{1-\alpha}} \right) + \frac{o(T_0-t)}{(T_0-t)^\alpha} W(x, t),  \quad (x, t) \in D \times [0, T_0)
\] for some stationary vector fields $V \in C^{2, \gamma}_{loc}$ and $W \in L^\infty \cap C^{2, 1}_{x, t, loc}$ and fixed $\alpha  <1$. If $W=0$, then $v$ is called a self-similar solution, abbreviated as SSS.
\end{defn}

\begin{rem}
One may extend the concept of ASSS by making a different scaling pattern for the angular component $V^\theta$ of the profile $V$ or making no extra assumption on the original angular component $v^\theta$ of the solution at all. For example, one can assume
$v^\theta =\frac{1}{(T_0-t)^\beta} V^\theta \left( \frac{x}{(T_0-t)^{1-\beta}} \right)$ with $0 \le \beta < \alpha$, with $\beta=0$ being the most natural choice for off axis considerations.  Since $v^\theta$ is bounded off the axis in our setting, as an error, after blowing up, it will disappear in any case.
\end{rem}

Throughout the paper, we make the following
basic assumptions for the solutions.

\noindent {\it Basic assumptions:}

{\it The vector field $v=v(x, t)$, $(x, t) \in D \times [0, T_0)$, $T_0
> 0$ is a smooth solution to the three-dimensional  axially symmetric Euler equations (ASE). Here $D$ is a finite cylinder around the $x_3$ axis and $v$ satisfies the no penetration boundary condition on the sides of $\partial D$, i.e. $v \cdot n=0$ where $n$ is the outward normal of $\partial D$, and $x_3$ periodic condition on the top and bottom of $\partial D$.}

 Now we are ready to state
the main result of the paper, together with its associated  corollary and  proposition.

\begin{thm}
\label{thsingstruc}
Let $v=v(x, t)$ be a smooth solution of ASE satisfying the Basic assumptions. Suppose $v$ develops a singularity at $t=T_0$ and off the $x_3$ axis. Then:

(a). for $\gamma \in (0, 1)$,  let $\{\td v_k\}$ be a blow-up sequence with locally bounded $C^{2, 1, \gamma}_{x, t}$ norms centered at anchor points of a singularity off the axis. Then any local $C^\gamma_{x, t}$ limit is a two dimensional, nontrivial  ancient solution of the Euler equation in either $\mathbb{R}^2 \times (-\infty, 0]$ or $\mathbb{R}^2_+ \times (-\infty, 0]$; if the anchor points are near maximal, then the limit solution is also bounded.

(b). in particular, for $\alpha  \in (0, 1)$, suppose $v$ is an asymptotic self-similar solution with a $C^{2, \gamma}$ profile and $C^{2, 1 \gamma}_{x, t}$ error, and  $|y| |\nabla V (y)|$ and $|y|^2 |\nabla^2 V (y)|$ are $C^\gamma$ functions. For $\gamma^- \in (0, \gamma)$,  any  blow-up sequence near maximal points off the  axis sub-converges in local $C^{\gamma^-}_{x, t}$ sense, to a  two dimensional, bounded, nontrivial  ancient solution to the Euler equation in either $\mathbb{R}^2 \times (-\infty, 0]$,  or in $\mathbb{R}^2_+ \times (-\infty, 0]$, which are functions of time only and decay to $0$ at negative infinite time.

\end{thm}

Under further assumptions such as choosing fixed spatial center point in scaling, we can allow the profile function to grow in sub-linear manner near infinity.

\begin{cor}
\lab{cornossb} Let $\alpha  \in (0, 1)$.
(1).
There does not exist asymptotic self-similar blow-up solutions with the following properties. (a) The profile $V$ and error are $C^{2, 1, \gamma}_{x, t}$  functions and $|y| |\nabla V (y)|$ and $|y|^2 |\nabla^2 V (y)|$ are $C^\gamma$ functions; (b).  near maximum points are within finite distance of the side boundary after blow-up; (c)  the vertical component of limiting profile is odd in $x^{(3)}$.

(2).
There does not exist asymptotic self-similar blow-up solutions with the following properties. (a). The spatial center of the blow-up sequence is at the side boundary point $p_0=(1, 0, 0)$ in the $(r, \theta, x^{(3)})$ coordinate.  (b). The profile $V$ and error are local $C^{2, 1, \gamma}_{x, t}$  functions; (c). anchor points are within finite distance of the side boundary after blow-up; (d). $|V(y)|$ is sub-linear near infinity and $|\nabla V|$ is bounded ; (e).  the vertical component of $V$ is odd in $x^{(3)}$.
\end{cor}

\begin{rem}
The growth conditions in part (2) of the Corollary are generous in that the profile of the velocity  is allowed to grow in sub-linear manner near infinity and the vorticity is allowed to be just bounded.

On the other hand , since there are numerous ancient and stationary solutions to the 2 dimensional Euler equation, it would be interesting to find one of them which is stable in certain weighted space. This would yield a velocity blow up result for the ASE.
\end{rem}

Next we turn to the  vorticity blow up scenario in \cite{LH2}, which is addressed in Propositions and remarks below, where the exponent $\alpha<0$.
The concept of ASSS in Definition \ref{defASSS} also covers the case $\alpha<0$ with suitable adjustment of the exponents on the components of the velocity. Some common ones proposed to show vorticity blow-up are: for $\alpha<0$ and $\beta>0$, $(x, t) \in D \times [0, T_0)$,
\be
\lab{LHsc}
\al
v=v(x, t)& = v^\theta e_\theta + v^r e_r + v^{(3)} e_3=\frac{1}{(T_0-t)^\alpha} \Theta \left( \frac{x}{(T_0-t)^{1-\alpha}} \right) e_{\theta} \\
&\qquad + \frac{(T_0-t)^\beta}{(T_0-t)^\alpha}  V^r\left( \frac{x}{(T_0-t)^{1-\alpha}} \right) e_r \\
&\qquad  + \frac{(T_0-t)^\beta}{(T_0-t)^\alpha}  V^{(3)} \left( \frac{x}{(T_0-t)^{1-\alpha}} \right) e_3 + \frac{o(T_0-t)}{(T_0-t)^\alpha} W(x, t).
\eal
\ee When $t \to T_0$, we notice that $v^r$ and $v^{(3)}$ become lower order terms since their vanishing order is $\beta$ more than that of $v^{\theta}$. Under this scenario, the vorticity would blow up at $T_0$ but the velocity stays bounded. Alternatively, one can also shift the factor $(T-t)^\beta$ around and assume
 \be
\lab{LHsc2}
\al
v=v(x, t)& = v^\theta e_\theta + v^r e_r + v^{(3)} e_3=\frac{1}{(T_0-t)^{\alpha+\beta}} \Theta \left( \frac{x}{(T_0-t)^{1-\alpha}} \right) e_{\theta} \\
&\qquad + \frac{1}{(T_0-t)^\alpha}  V^r\left( \frac{x}{(T_0-t)^{1-\alpha}} \right) e_r \\
&\qquad  + \frac{1}{(T_0-t)^\alpha}  V^{(3)} \left( \frac{x}{(T_0-t)^{1-\alpha}} \right) e_3 + \frac{o(T_0-t)}{(T_0-t)^\alpha} W(x, t).
\eal
\ee

 The following proposition appears to show that most of this kind of vorticity blow up does not occur either. The exceptional case is when $2 \beta=1-\alpha$.

  The main point is that if the three components of the velocity vanish at different speed at $T_0$, then a suitable scaling limit of the solutions will be simple enough to allow us to reach definite conclusion.

\begin{proposition}
\lab{prLH}
(i) There does not exist asymptotic self-similar blow-up solutions \eqref{LHsc} with the following properties. (a). The spatial center of the blow-up sequence is at the side boundary point $p_0=(1, 0, 0)$ in the $(r, \theta, x^{(3)})$ coordinate.  (b). The profiles $\Theta$, $V^r$, $V^{(3)}$  are nontrivial, local $C^{2, \gamma}_{x}$ functions and  the error $W$ are local $C^{2, 2, \gamma}_{x, t}$  functions.

(ii)  There does not exist asymptotic self-similar blow-up solutions \eqref{LHsc2} with the following properties. (a). The spatial center of the blow-up sequence is at the side boundary point $p_0=(1, 0, 0)$ in the $(r, \theta, x^{(3)})$ coordinate.  (b). The profiles $\Theta$, $V^r$, $V^{(3)}$  are nontrivial, local $C^{2, \gamma}_{x}$ functions and  the error $W$ are local $C^{2, 2, \gamma}_{x, t}$  functions. (c). $v^\theta$ is odd in $x_3$; (d) $-2(\beta/\alpha)+(1/\alpha)-1>0$.

(iii) Suppose there exists asymptotic self-similar blow-up solutions \eqref{LHsc2} with
\be
\lab{crtab} -(2\beta/\alpha)+(1/\alpha)-1=0, \quad i.e. \quad 2 \beta = 1-\alpha.
\ee Then a suitable blow up converges in local $C^\gamma_{x, t}$  sense to an exact self similar solution of the two dimensional Boussinesq's equations  in the half plane.
\end{proposition}

  Note this proposition does not exclude a  proposed self similar blow up scenario \cite{LH1} when $2 \beta = 1-\alpha$, which is  expected to be true in the literature with $v^\theta, \omega^\theta$ being odd in $x^{(3)}$ variable. For example, the expected values are $\alpha \approx -2$ and $\beta \approx 1.5$ satisfying \eqref{crtab}  as stated in the recent \cite{CH1} and \cite{WLG-SB}.  This situation will be discussed in Propositions
  \ref{prLH00}, \ref{prLH0}  and the remarks in Section 2. The case when $-2\beta/\alpha+(1/\alpha)-1<0$ is similar to Corollary \ref{cornossb} part (2). See Proposition \ref{prLHextra} below.

Some relaxations of the conditions in the proposition (ii) are also possible. See Remark
\ref{rmk2} below.

Here are some notations to be used frequently. We use
$v=v(x, t)$ to denote a solution (velocity field) to the ASE. Here $(x, t)$ is
a point in space-time. Unless stated otherwise,  we use $r, r_0, r_k$ to denote the
distance between points $x, x_0, x_k$ in space and the $x_3$ axis (rotational axis)
respectively. $SSS$ stands for self similar solutions. If $x=(x^{(1)}, x^{(2)}, x^{(3)})$ then $x'$ denotes $(x^{(1)}, x^{(2)}, 0)$. The theorem, corollary and proposition will be proven in Section 2.

\section{Proof of Theorem \ref{thsingstruc} and the propositions}

This section is divided into two parts. In Part I,  we prove Theorem \ref{thsingstruc} and in part II, we prove Proposition \ref{prLH} and state and prove three new propositions concerning vorticity blow up.
\medskip

{\it Part I.}

\proof (of Theorem \ref{thsingstruc}).

It is well-known that local in time existence and uniqueness of smooth solutions of \eqref{EL} in smooth bounded domains with no penetration boundary condition have been proven in \cite{EM}, \cite{BB} and
\cite{Te}.

Let us prove part (a) first.

We suppose the singularity happens at $t=1$ the first time. Let $(x_k, t_k)$ be a sequence of anchor  points and
\begin{equation}
\lab{defqk}
Q_k=|v(x_k, t_k)|.
\end{equation} Here the coefficient is taken as $1$ for simplicity. Since we assume the singularity occurs away from the axis or at the side boundary of the cylinder $D$, we can assume, without loss of generality, that
\be
\lab{xk'>1/2}
|x_k'| \ge 1/2.
\ee

Fix a number $\alpha \in (0, 1)$, define the scaled function
\begin{equation}\label{vktild}
\tilde v_k= \tilde v_k(\tilde x, \tilde t) = Q_k^{-1} v(Q_k^{-(1-\alpha)/\alpha} \tilde x + x_k,
Q_k^{-1/\alpha} \tilde t + t_k), \quad \tilde t \le 0.
\end{equation}
Then $\tilde v_k$ is a solution of the Euler equation in
the slab $D_k \times [ - Q_k^{1/\alpha}/2, 0]$. Here
\be
\lab{defdk}
D_k = \{\tilde x \, | \,  Q_k^{-(1-\alpha)/\alpha} \tilde x + x_k \in D \}.
\ee As $Q_k \to \infty$, we will see that $D_k$ will expand in size and eventually either become the whole space or the half space for the variable $\tilde x$.
Moreover,
by the definition of anchor points in Definition \ref{defbaopuops} on $Q_k$, we know that
\be
\lab{vkt<4}
|\tilde v_k(\tilde x, \tilde t)| \le \lambda(|\td{x}|^2 +|\td t|)  |\tilde v_k(0, 0)|
=\lambda(|\td{x}|^2 +|\td t|)
\ee
whenever defined. So the sequence of blow-up solutions $\tilde v_k$ are uniformly locally bounded.

In the standard basis for $\mathbb{R}^3$, let $x_k=(x_{k}^{(1)},
x_{k}^{(2)}, x_{k}^{(3)})$ with the third component being the one for
the vertical axis, and let $\xi_k = (0, 0, x^{(3)}_{k})$. Since the sequence of
vectors
\[
\{(x_k-\xi_k)/|x_k-\xi_k|\}
\]  are unit ones, there exists a
subsequence, still labeled by $k$, which converges to a unit
vector $\zeta=(\zeta_1, \zeta_2, 0)$.  We use the three vectors
\begin{equation}\nonumber
\zeta=(\zeta_1, \zeta_2, 0), \,  \zeta^\bot=(-\zeta_2, \zeta_1, 0), \,  e_3=(0, 0, 1)
\end{equation}
as
the basis of a new coordinate.  Since this basis is obtained by a
rotation around the vertical axis, we know $v$ is invariant.  From now
on, when we mention the coordinates of a point $x=(x^{(1)}, x^{(2)}, x^{(3)})$, we mean to use the
new basis with the same origin, i.e. $x=x^{(1)} \zeta + x^{(2)} \zeta^\bot +x^{(3)} e_3$. We still use $(r, \theta,  x^{(3)})$ to
denote the variables in the cylindrical system corresponding to the point $x$ in
this new basis, namely $r= \sqrt{(x^{(1)})^2 + (x^{(2)})^2}, \, \theta= \tan^{-1} (x^{(2)}/x^{(1)})$.

Let $\beta_k$ be a sequence of positive numbers converging to $0$ as $k \to \infty$ at a rate which is slower than that of $Q^{-(1-\af)/\af}_k$. For $x \in B(x_k, \beta_k) \cap D$, we recall that $\theta$ is the longitude angle
between $x'$ and $\zeta$. Then
\begin{equation}\label{thetakto0}
\cos \theta =
 \frac{(x-\xi_k) \cdot
\zeta}{|x - (0, 0, x_3)|} = \frac{(x_k-\xi_k) \cdot \zeta}{|x_k
-\xi_k|} + \frac{O(\beta_k)}{r_k} \to 1, \qquad k \to \infty.
\end{equation}

For the solution of the Euler equation $v=v(x, t)$  in $[B(x_k, \beta_k) \cap D] \times [t_k-\beta^2_k,
t_k]$, recall from \eqref{vktild} the re-scaled solution
\begin{equation}\nonumber
\tilde v_k= \tilde v_k(\tilde x, \tilde t) = Q_k^{-1} v(Q_k^{-(1-\alpha)/\alpha} \tilde x + x_k,
Q_k^{-1/\alpha} \tilde t + t_k)
\end{equation}
where
$x=Q_k^{-(1-\af)/\af } \tilde x + x_k$ and $t=Q_k^{-1/\af} \tilde t + t_k$. Then
for $x=(x^{(1)}, x^{(2)}, x^{(3)})$ and $\tilde x=(\tilde x^{(1)},
\tilde x^{(2)}, \tilde x^{(3)})$, $r= \sqrt{(x^{(1)})^2 + (x^{(2)})^2 }$,  we have
\begin{equation}\label{vkdaoshu}
\begin{cases}
\partial_r v(x, t) = \partial_{x^{(1)}} v(x, t)
  \cos \theta + \partial_{x^{(2)}} v(x, t)  \sin \theta\\
\quad\quad = Q^{1/\af}_k \partial_{\tilde x^{(1)}} \tilde v_k(\tilde x,
  \tilde t) \cos \theta + Q^{1/\af}_k \partial_{\tilde x^{(2)}} \tilde v_k (\tilde x, \tilde t) \sin
  \theta\\[-4mm]\\
\partial_{x^{(3)}} v(x, t) = Q^{1/\af}_k \partial_{\tilde x^{(3)}} \tilde v_k (\tilde x, \tilde t) \\
\partial_t v(x, t) = Q^{1+1/\af}_k \partial_{\tilde t} \tilde v_k(\tilde x, \tilde t).
\end{cases}
\end{equation}

For the pressure $p=p(x, t)$, the scaled ones are
\begin{equation}\nonumber
\tilde p_k =
 \tilde p_k(\tilde x, \tilde t) = Q_k^{-2} p(Q_k^{-(1-\alpha)/\alpha} \tilde x + x_k,
Q_k^{-1/\alpha} \tilde t + t_k).
\end{equation}
Therefore
\begin{equation}\label{pdaoshu}
\partial_r p(x, t) =Q^{1+1/\af}_k \partial_{\tilde x^{(1)}} \tilde p_k(\tilde x, \tilde t)  \cos \theta +
Q^{1+1/\af}_k \partial_{\tilde x^{(2)}} \tilde p_k (\tilde x, \tilde t)
\sin \theta
\end{equation}

Writing $v= v^r e_r +  v^\theta e_\theta + v^{(3)} e_3$, then
\be\label{driftdaoshu}
\al
&v^r \partial_r v^r + v^{(3)} \partial_{x^{(3)}} v^r \\
&= Q^{1+1/\af}_k \bigg[
 \tilde v^r_k(\tilde x, \tilde t) \partial_{\tilde x^{(1)}} \tilde
  v^r_k(\tilde x, \tilde t)  \cos \theta \\
 &\qquad +\ \tilde v^r_k(\tilde x, \tilde t)
  \partial_{\tilde x^{(2)}} \tilde v^r_k (\tilde x, \tilde t) \sin
  \theta + \tilde v^{(3)}_k \partial_{\tilde x^{(3)}} \tilde
  v^r_k(\tilde x, \tilde t) \bigg].
  \eal
\ee
We substitute the above identities into the equation for $v^r$ in \eqref{axi-EL}:
\begin{equation}\nonumber
-(b\cdot\nabla)v^r+
\frac{(v^{\theta})^2}{r}-\frac{\partial p}{\partial
r}-\frac{\partial v^r}{\partial t}=0,
\end{equation}
then we arrive at
\be
\lab{eqtvro}
\al
-( \tilde v^r_k \partial_{\tilde x^{(1)}} + \tilde
  v^{(3)}_k \partial_{\tilde x^{(3)}} ) \tilde v^r_k -\partial_{\tilde
  x^{(1)} }\tilde p_k - \partial_{\tilde t} \tilde v^r_k
+ \frac{( \tilde v^\theta_k)^2}{Q_k^{[(1/\alpha)-1]} r}+ O_r(\theta)=0.
\eal
\ee
Here
\[
\al
O_r(\theta)&= \tilde v^r_k(\tilde x, \tilde t) \partial_{\tilde x^{(1)}} \tilde
  v^r_k(\tilde x, \tilde t) (1- \cos \theta)  -\ \tilde v^r_k(\tilde x, \tilde t)
  \partial_{\tilde x^{(2)}} \tilde v^r_k (\tilde x, \tilde t) \sin
  \theta\\
  &\qquad + \partial_{\tilde x^{(1)}} \tilde p_k(\tilde x, \tilde t) (1- \cos \theta) -
 \partial_{\tilde x^{(2)}} \tilde p_k (\tilde x, \tilde t)
\sin \theta
\eal
\] which only contains terms that
vanish when $\theta \to 0$ as $k \to \infty$.  In particular all
terms involving the derivative with respect to $\tilde x^{(2)}$
are included in $O_r(\theta)$.
 Similarly
\be
\lab{eqtv3o}
-( \tilde v^{r}_k \partial_{\tilde x^{(1)}} + \tilde
v^{(3)}_k \partial_{\tilde x^{(3)}} ) \tilde v^{(3)}_k
-\partial_{\tilde x^{(3)} }\tilde p_k -
\partial_{\tilde t} \tilde v^{(3)}_k  + O_3(\theta)=0
\ee where
\[
 O_3(\theta)=\tilde v^r_k(\tilde x, \tilde t) \partial_{\tilde x^{(1)}} \tilde
  v^{(3)}_k(\tilde x, \tilde t) (1- \cos \theta)  -\ \tilde v^{(3)}_k(\tilde x, \tilde t)
  \partial_{\tilde x^{(2)}} \tilde v^{r}_k (\tilde x, \tilde t) \sin
  \theta.
\]In the above we have employed the fact that the $C^{0, 0, \gamma}_{x, t}$ norm of $\nabla \tilde p_k$ are bounded on compact sets thanks to the Euler equations:
\[
\nabla \tilde p_k = - \partial_t \tilde v_k - \tilde v_k \cdot \nabla \tilde v_k.
\]

Notice that for points $x=(x^{(1)}, x^{(2)}, x^{(3)}) \in B(x_k, \beta_k) \cap D$, we have $r=|x'| \ge 1/4$ since $|x'_k| \ge 1/2$ and $\beta_k \to 0$. Therefore, due to $\alpha \in (0, 1)$, we have $Q_k^{[(1/\alpha)-1]} r \to \infty$ and $\theta = \tan^{-1}(x^{(2)}/x^{(1)}) \to 0$ as $k \to \infty$.

 Letting $k \to \infty$ and remembering  that $|\tilde v_k(\cdot, \cdot)| \le \lambda(\cdot)$ from \eqref{vkt<4}, we know by assumption  that  $\{ \tilde v_k \}$ sub-converges in local $C^{\gamma}_{x, t}$ sense to the bounded vector fields
\be
\lab{vktlim}
\tilde v=\tilde v(\tilde x, \tilde t) = \tilde v^{(1)} \zeta + \tilde v^{(2)} \zeta^\bot + \tilde v^{(3)} e_3,
\ee and $\nabla \tilde p_k$ sub-converges in $C^{\gamma}_{x, t}$ sense to a vector field $\nabla \tilde p$.
In addition, the functions $\tilde v^{(1)}$ and $\tilde v^{(3)}$, as bounded, local $C^{2, 1, \gamma}_{x, t}$ functions,  satisfy  the equations
\be
\lab{eqvt2d}
\begin{cases}
 ( \tilde v^{(1)} \partial_{\tilde x^{(1)}} + \tilde
v^{(3)} \partial_{\tilde x^{(3)}} ) \tilde v^{(1)}
+\partial_{\tilde x^{(1)} }\tilde p +
\partial_{\tilde t} \tilde v^{(1)} =0,\\
 ( \tilde v^{(1)} \partial_{\tilde x^{(1)}} + \tilde
v^{(3)} \partial_{\tilde x^{(3)}} ) \tilde v^{(3)}
+\partial_{\tilde x^{(3)} }\tilde p +
\partial_{\tilde t} \tilde v^{(3)} =0.
\end{cases}
\ee Observe that the domain for $\tilde v_k$ contains the space time region
\[
\Omega_k \equiv [B(0, \beta_k Q^{(1-\af)/\af}_k) \cap D_k] \times [-\beta^2_k Q^{1/\af}_k/2, 0].
\]For the $\tilde x$ variable, $D_k$ is a cylinder with scale $Q^{(1-\af)/\af}_k$ and  reference point $0$ which corresponds to the point $x_k$ in the un-scaled coordinates. The distance from this reference point $0$ to the side boundary of $D_k$ is $Q^{(1-\af)/\af}_k (1-|x'_k|)$. We have two cases to consider. Case one is that the sequence  $\{ Q^{(1-\af)/\af}_k (1-|x'_k|) \}$ stays bounded. Then, since
 $\beta_k Q^{(1-\af)/\af}_k \to \infty$ by our choice, the domains $\Omega_k$ expands to $\mathbb{R}^3_+ \times (-\infty, 0]$. Here $\mathbb{R}^3_+$ is the half space $\tilde x^{(1)} \le c$ for some constant $c \ge 0$ which will be made $0$ after a translation. Case 2 is that a subsequence of $\{ Q^{(1-\af)/\af}_k (1-|x'_k|) \}$ diverges to positive infinity. Then it is clear that $\Omega_k$ expands to $\mathbb{R}^3 \times (-\infty, 0]$.
 Hence the domain for the variables in equation \eqref{eqvt2d} is either $\mathbb{R}^3_+ \times (-\infty, 0]$ or $\mathbb{R}^3 \times (-\infty, 0]$.

Next, we prove that $\tilde v^{(2)} \equiv 0$.
First, we recall the fact that the component $v^\theta$  in the original solution $v$ is bounded when the points are bounded away from the $x_3$ axis. Indeed, from \eqref{axi-EL}, for $t<1$, the function
\[
\Gamma \equiv r v^\theta (x, t)
\]satisfies, in the classical sense, the equation
\be
\lab{eqgam}
\partial_t \Gamma + b \nabla \Gamma =0.
\ee Since $D$ is a bounded domain, the no penetration boundary condition and integration by parts infer
\[
\int_D |\Gamma(x, t)|^{2n} dx = \int_D |\Gamma(x, 0)|^{2n} dx, \quad n=1, 2, 3, ....
\]This implies
\[
\Vert \Gamma(\cdot, t) \Vert_\infty = \Vert \Gamma(\cdot, 0) \Vert_\infty.
\]Hence
\be
\lab{vthjie}
|v^\theta(x, t)| \le \frac{\Vert \Gamma(\cdot, 0) \Vert_\infty}{r}, \quad t<1.
\ee From this and the relation
 \[
 \tilde v^\theta_k = Q^{-1}_k v^\theta(Q_k^{-(1-\alpha)/\alpha} \tilde x + x_k,
Q_k^{-1/\alpha} \tilde t + t_k),
 \] we know that for $(\tilde x, \tilde t) \in \Omega_k$, the following bounds hold
 \[
 |\tilde v^\theta_k | \le C Q^{-1}_k.
 \]Here we have used the fact that for points $(\tilde x, \tilde t) \in \Omega_k$,  the corresponding original spatial variable $x$  is bounded away for the vertical axis.   This implies that $\tilde v^{(2)}$, as the $L^\infty$ limit of $\tilde v^\theta_k$ is zero.

Finally, we need to show that $\tilde v^{(1)}$ and $\tilde v^{(3)}$
are independent of the variable $\tilde x^{(2)}$.  To prove it,
let us observe that by axial symmetry,  $\partial_\theta \td v^r_k=\partial_\theta
\td v^{(3)}_k=0$. Hence
\begin{equation}\nonumber
-\partial_{x^{(1)}}\td v_k^r \sin \theta + \partial_{x^{(2)}}\td v_k^r
\cos \theta = -\partial_{x^{(1)}}\td v_k^{(3)} \sin \theta +
\partial_{x^{(2)}}\td v_k^{(3)} \cos \theta = 0.
\end{equation}
This implies, in the classical sense, that
\begin{equation}\nonumber
\partial_{\tilde x^{(2)}}  \tilde v_k = \partial_{\tilde x^{(1)} } \tilde v_k \tan
\theta. \end{equation}
Taking $k \to \infty (\theta \to 0)$ we see
the desired result is true. Notice that the axially symmetric divergence free condition also becomes the 2-dimensional divergence free condition after taking the limit.

Hence we have proven part (a) of the theorem, i.e., the local $C^{\gamma}_{x, t}$ blow-up limit of potential singularity near maximal points off the vertical axis is a bounded, two-dimensional ancient solution of the Euler equation in either $\mathbb{R}^2 \times (-\infty, 0]$ or $\mathbb{R}^2_+ \times (-\infty, 0]$. Note that we have changed $\mathbb{R}^3$ to $\mathbb{R}^2$ due to the absence of the variable $\tilde x^{(2)}$. This concludes the proof of part (a) of the theorem.

In summary the limit of blow-up solutions $\tilde v=\tilde v(\tilde x^{(1)}, \tilde x^{(3)}) = \tilde v^{(1)} \zeta + \tilde v^{(3)} e_3 \in C^{2, 1, \gamma}_{x, t}$ satisfies the two dimensional Euler equation
\be
\lab{eqvtlim2d}
\begin{cases}
	( \tilde v^{(1)} \partial_{\tilde x^{(1)}} + \tilde
	v^{(3)} \partial_{\tilde x^{(3)}} ) \tilde v^{(1)}
	+\partial_{\tilde x^{(1)} }\tilde p +
	\partial_{\tilde t} \tilde  v^{(1)} =0,\\
	( \tilde v^{(1)} \partial_{\tilde x^{(1)}} + \tilde
	v^{(3)} \partial_{\tilde x^{(3)}} ) \tilde v^{(3)}
	+\partial_{\tilde x^{(3)} }\tilde p +
	\partial_{\tilde t}\tilde  v^{(3)} =0
\end{cases}
\ee in $\tilde D \times (-\infty, 0]$ where $\tilde D = \mathbb{R}^2$ or $ \mathbb{R}^2_+$ with no penetration boundary condition. Moreover $\tilde v$ is not identically zero since $|\tilde v(0, 0)|=1$.

 \medskip

Now we prove part (b).

In part (a) we assumed the blow-up limit converges  in local $C^\gamma_{x, t}$ sense.  In this part (b), we assume that the original solution is an asymptotic self-similar one whose profile and error are bounded $C^{2, \gamma}$ function. Then we actually can prove convergence  in the local $C^{\gamma^-}_{x, t}$ sense. Here $\gamma^- \in (0, \gamma)$. In addition,  for the limit 2-dimensional ancient solutions,  the no penetration boundary condition (if the limiting spatial domain is $\mathbb{R}^2_+$), are also preserved and it decays to $0$ at $\tilde t = -\infty$.

The rest of the proof is divided into two steps.

{\it Step 1.} For clarity of presentation we first assume the error term is zero.

In this case, for some $ \alpha \in (0, 1)$, the solution is given by
\be
\lab{defvss}
v=v(x, t) = \frac{1}{(1-t)^\alpha} V \left( \frac{x}{(1-t)^{1-\alpha}} \right),  \quad (x, t) \in D \times [0, 1),
\ee where $V$ is a  bounded, $C^{2, \gamma}$ vector field with 3 spatial variables by assumption. Here for simplicity, we assume the first singular time is $t=1$. We suppose, the points $(x_k, t_k) \in D \times [0, 1)$ are near maximum points for $|v|$ in $D \times [0, t_k]$ and $t_k \uparrow 1$. Since we are only concerned with off-axis blow-up, we can assume without loss of generality that $|x'_k| \ge 1/2$. Since $|V|$ is bounded, we know that $|v(x_k, t_k)|$ is comparable to $(1-t_k)^{-\alpha}$.

Write
\[
Q_k=(1-t_k)^{-\alpha}.
\]As in part (a), we consider the scaled solution
\begin{equation}\label{vktild2}
\tilde v_k= \tilde v_k(\tilde x, \tilde t) = Q_k^{-1} v(Q_k^{-(1-\alpha)/\alpha} \tilde x + x_k,
Q_k^{-1/\alpha} \tilde t + t_k), \quad \tilde t \le 0.
\end{equation} Due to \eqref{defvss}, $\tilde v_k$ takes the following special form
\begin{equation}
\label{vktild3}
\al
\tilde v_k &= \frac{Q_k^{-1}}{(\, 1-t_k-Q_k^{-1/\alpha} \,  \tilde t \, )^{\alpha} } V\left(\frac{Q_k^{-(1-\alpha)/\alpha} \tilde x + x_k}{[\, 1-t_k-Q_k^{-1/\alpha} \,  \tilde t \, ]^{1-\alpha} } \right)\\
&=\frac{1}{(\, 1 - \tilde t \, )^{\alpha} } V\left(\frac{\tilde x + x_k Q_k^{(1-\alpha)/\alpha}}{(1-\tilde t)^{1-\alpha}}  \right)\\
&=\frac{1}{(\, 1 - \tilde t \, )^{\alpha} } V\left(\frac{\tilde x + z_k}{(1-\tilde t)^{1-\alpha}}  \right).
\eal
\end{equation} Here for simplicity we have written
$z_k= x_k Q_k^{(1-\alpha)/\alpha}$.

 We will prove that $\tilde v_k$ converges in local $C^{\gamma^-}_{x, t}$ sense.

From \eqref{vktild3} and our assumption on $V$, we know that there is a uniform constant such that
\be
\lab{vkc2gam}
\Vert \tilde v_k(\cdot, \tilde t) \Vert_{C^{2, \gamma}} \le C, \quad \forall \tilde t \le 0.
\ee Since $\tilde v_k$ is a smooth solution of the Euler equation in the scaled up cylinder $D_k$, by differentiating \eqref{vktild3}, we have
\be
\lab{ptvkgam}
\al
&\Vert \partial_{\td t} \tilde  v_k(\cdot, \tilde t) \Vert_{C^{0}} \\
&\le \frac{1-\af}{(\, 1 - \tilde t \, )^2 } |\td x + z_k| \, | \nabla V| \left(\frac{\tilde x + z_k}{(1-\tilde t)^{1-\alpha}}  \right) + \frac{\af}{(\, 1 - \tilde t \, )^{1+\alpha} } |V|\left(\frac{\tilde x + z_k}{(1-\tilde t)^{1-\alpha}}  \right) \\
&\le \frac{C}{(\, 1 - \tilde t \, )^{1+\af}} , \quad \forall \tilde t \le 0.
\eal
\ee Here we have used the assumption that $|y| |\nabla V(y)|$ is bounded.  Similarly, using also the boundedness assumption of $|y|^2 |\nabla^2 V(y)|$, we see that
\[
\Vert \partial_{\td t} \nabla \tilde  v_k(\cdot, \tilde t) \Vert_{C^{0}} \le C, \quad \forall \td t \le 0;
\]
\[
\Vert \partial_{\td t} \partial_{\td t}  \tilde  v_k(\cdot, \tilde t) \Vert_{C^{0}} \le C, \quad \forall \td t \le 0;
\]
\[
\Vert    \nabla \tilde  p_k(\cdot, \tilde t) \Vert_{C^{0}} \le C, \quad \forall \td t \le 0.
\]The third bound is due to the Euler equation $\nabla \tilde p_k= - \tilde v_k \cdot \nabla \tilde v_k - \partial_{\td t} \tilde v_k$ and \eqref{vkc2gam}, \eqref{ptvkgam}.

Differentiating the Euler equation
\be
\lab{euler2}
 \partial_{\td t} \tilde v_k = - \nabla \tilde p_k  - \tilde v_k \cdot \nabla \tilde v_k
\ee
with respect to time, we find
\be
\lab{euler2t}
\partial_{\td t} \partial_{\td t} \tilde v_k =- \nabla \partial_{\td t} \tilde p_k  - \tilde v_k \cdot \nabla \partial_{\td t} \tilde  v_k - \partial_{\td t} \tilde  v_k \cdot \nabla  \tilde  v_k.
\ee Hence
\[
\Vert  \nabla \partial_{\td t} \td p_k \Vert_{C^{0}(D_k)} \le C_5,
\]

Putting together, we deduce
\be
\lab{vktt}
\Vert  \partial_{\td t} \tilde v_k \Vert_{C^{0, \gamma}(D_k)}+ \Vert \partial_{\td t} \nabla \tilde  v_k(\cdot, \tilde t) \Vert_{C^{0}} +
 \Vert  \nabla \partial_{\td t} \td p_k \Vert_{C^{0}(D_k)}+\Vert \partial_{\td t} \partial_{\td t} \tilde v_k \Vert_{C^{0}(D_k)} \le C_6
\ee for a uniform constant $C_6$.

 Bound \eqref{vktt} gives us enough regularity in time to take limit in local $C^{\gamma^-}_{x, t}$ sense. Notice that the $C^{2, \gamma}$ spatial regularity is guarantied by our assumption on the profile function $V$.
 Since $\tilde v_k$ is a solution of the Euler equation, by Ascoli-Arzela theorem, we conclude that its local $C^{\gamma^-}_{x, t}$  limit vector field
\be
\lab{limtldv}
\tilde v \equiv \frac{1}{(\, 1 - \tilde t \, )^{\alpha} } \tilde V\left(\frac{\tilde x }{(1-\tilde t)^{1-\alpha}}, \tilde  t \right)
\ee is also a solution of the Euler equation.
Indeed, let $t_i$ be an enumeration of the rational numbers in $[-1, 0]$. For each $i=1, 2, 3, ...$
consider the translated stationary vector fields
\[
V_{k, i}=V_{k, i}\left(\frac{\tilde x}{(1-\tilde t_i)^{1-\alpha}} \right) \equiv V\left(\frac{\tilde x + z_k}{(1-\tilde t_i)^{1-\alpha}} \right).
\]By \eqref{vkc2gam}, for $i=1$, there is a sequence of integers $\{1_k\}$ such that $V_{1_k, 1}$ converges, in local $C^{2, \gamma^-}$ sense, to a vector field $\tilde V_1 = \tilde V_1\left(\frac{\tilde x}{(1-\tilde t_1)^{1-\alpha}} \right)$. By \eqref{vkc2gam} again, for $i=2$, there is a subsequence of integers $\{2_k\} \subset \{1_k\}$ such that $V_{2_k, 2}$ converges, in local $C^{2, \gamma^-}$ sense, to a vector field $\tilde V_2 = \tilde V_2\left(\frac{\tilde x}{(1-\tilde t_2)^{1-\alpha}} \right)$. Repeating this process and choosing the diagonal sequence, denoted by $\{ V_{k_k}\}$, we see that
$V_{k_k}$ converges, in local $C^{2, \gamma^-}$ sense, to a vector field $\tilde V = \tilde V\left(\frac{\tilde x}{(1-\tilde t_i)^{1-\alpha}}, t_i \right)$ for each $i=1, 2, 3, ...$. Using \eqref{ptvkgam}, we can extend the vector fields $\tilde V_i$ to all real numbers in $[-1, 0]$. The extended vector field is denoted by
\[
\tilde V= \tilde V\left(\frac{\tilde x}{(1-\tilde t)^{1-\alpha}}, \tilde t \right).
\]Then we can deduce, via \eqref{vktt}, as in the proof of the Ascoli-Arzela theorem , that the sequence
\[
\frac{1}{(1-\tilde t)^\alpha} V_{k_k}=\frac{1}{(1-\tilde t)^\alpha } V_{k_k}\left(\frac{\tilde x}{(1-\tilde t)^{1-\alpha}} \right)
\]converges, in local $C^{\gamma^-}_{x, t}$ sense, to $\frac{1}{(1-\tilde t)^\alpha} \tilde V$, $\tilde t \in [-1, 0]$.  Repeating this process, we obtain the convergence for all $\tilde t \le 0$.

According to part (a) of the theorem, $\tilde v$ in \eqref{limtldv} is an ancient solution of the 2-dimensional Euler equations on the plane or the half plane. Due to the $C^{\gamma^-}_{x, t}$ convergence, the no penetration boundary condition is also preserved in the half plane case.

 So the blow-up limit is a nontrivial ancient solution $\tilde v$ of the Euler equation in $\tilde D= \mathbb{R}^2$ or $\mathbb{R}^2_{+}$, which are in the space $C^{2, 1, \gamma}_{x, t}$. We mention that for each fixed $\td t$, $\td v_k(\cdot, \td t)$ has uniformly bounded $C^{2, \gamma}$ norm by assumption. Therefore, even though the convergence is in the $C^{\gamma^-}_{x, t}$ sense  (including $\nabla \td v_k$), the limit function $\td v(\cdot, \td t)$ is still in $C^{2, \gamma}$ after extracting a local $C^{2, \gamma}$ convergent subsequence for the fixed time slice.
 In the case there is boundary, $\tilde v$ also satisfies the no penetration boundary condition.
Let $\omega$ be the (scalar) vorticity of $\tilde v$, which satisfies the transport equation
\be
\lab{2dvtx}
\partial_{\tilde t} \omega + \tilde v \nabla \omega=0.
\ee Due to the boundary condition, the flow lines of $\tilde v$ do not cross the boundary of $\mathbb{R}^2_+$. Hence the sup norm of $\omega(\cdot, \tilde t)$ does not increase in time. But it is clear from \eqref{limtldv} that $\lim_{\tilde t \to -\infty} \Vert \omega(\cdot, \tilde t)\Vert_\infty=0$. Hence $\omega \equiv 0$. Therefore $\tilde v(\cdot, \tilde t)$ is a bounded harmonic function since
\[
\Delta \tilde v = -\nabla \times (\omega \, e^{(3)}) =0.
\]In case the domain $\tilde D$ is the full plane $\mathbb{R}^2$, then $v=v(t)$ by the Liouville theorem.

In case the domain $\tilde D$ is the half plane $\mathbb{R}^2_+$, say $\tilde x^{(1)} \le 0$, for simplicity we drop the tilde symbol for all variables for the rest of the proof. Also, in \eqref{eqvtlim2d},  we use $x_1$ to replace $\tilde x^{(1)}$, $x_2$ to replace $\tilde x^{(3)}$, $v_1$ to replace $\tilde v^{(1)}$ and $v_2$ to replace $\tilde v^{(3)}$. Recall that $\tilde v^{(2)}=0$.  Since $v_1=0$ at the boundary, by the Liouville theorem again $v_1=0$. The divergence free condition $\partial_1 v_1 + \partial_2 v_2=0$ implies that $v_2=v_2(x_1, t)$.
Using $0= \omega = \partial_2 v_1 - \partial_1 v_2$, we see that $\partial_1 v_2=0$ and hence $v_2=v_2(t)= c/(1-t)^\alpha$. Therefore
\[
v=(v_1, v_2)=(0,  c/(1-t)^\alpha).
\]  This proves part (b) of the theorem in case the error $W=0$.
\medskip

{\it Step 2.}  The error $W$ is not zero.

Since the error term tends to $0$ point-wise after blow-up, we can modify the proof in Step 1 without much difficulty.  The reason is that, in the presence of error,  we can decompose the scaled solutions in the sequence as the purely self-similar part $S_k$ plus the terms involving the error $E_k$:
\[
\al
\tilde v_k &= S_k + E_k \\
&\equiv \frac{1}{(\, 1 - \tilde t \, )^{\alpha} } V\left(\frac{\tilde x + z_k}{(1-\tilde t)^{1-\alpha}}  \right)+ \frac{o(1)}{(\, 1 - \tilde t \, )^{\alpha} } W\left( Q_k^{-(1-\alpha)/\alpha} \tilde x + x_k,
Q_k^{-1/\alpha} \tilde t + t_k \right).
\eal
\] In the scaled variables, the Euler equation becomes
\be
\lab{eeSkEk}
S_k \nabla S_k + \partial_{\tilde t} S_k + \nabla p_k +  S_k \nabla E_k + E_k \nabla S_k
 + E_k \nabla E_k + \partial_{\tilde t} E_k =0.
\ee

  We  choose a smooth, divergence free, compactly supported  test vector fields $\zeta_k=\zeta_k(\tilde x, \tilde t)$ in such a way that their supports are within a compact domain containing a fixed reference point and of fixed diameter in the $(\tilde x, \tilde t)$ space time and that their $C^{2, 1, \gamma}_{x, t}$ norms are uniformly bounded. Using $\zeta_k$ as a test function in \eqref{eeSkEk}, we see the pressure terms disappear, giving us
\be
\lab{skekint}
\int \int \left[ S_k \nabla S_k + \partial_{\tilde t} S_k \right] \zeta_x d\tilde x d\tilde t + \int \int \left[ S_k \nabla E_k + E_k \nabla E_k + E_k \nabla S_k + \partial_{\tilde t} E_k \right] \zeta_k d\tilde x d\tilde t=0.
\ee
 In \eqref{skekint}, the purely self-similar part $S_k$ was dealt with in Step 1, which sub-converges in $C^{\gamma^-}_{x, t}$ sense to $\tilde v$. The terms involving the error $E_k$ converges to $0$ since each one contains $o(1)$. We can choose $\zeta_k$ suitably so that they converge in $C^{2, 1, \gamma^-}_{x, t}$ sense, to a given test vector field $\zeta$ which is also compactly supported. After taking $k \to \infty$, we deduce from \eqref{skekint} that
 \[
 \int \int \left[ \tilde v \nabla \tilde v + \partial_{\tilde t} \tilde v \right] \zeta d\tilde x d\tilde t =0.
 \]
Therefore the purely self-similar part sub-converges in the $C^{\gamma^-}_{x, t}$  sense to a weak solution of the Euler equation $\tilde v$. Since this limit $\tilde v$ itself is a $C^{2, 1, \gamma}_{x, t}$ function, it must also be a classical solution. Then we can just proceed as in the last two paragraph of Step 1. This finishes the proof of the theorem.
\qed

\medskip

Next we give a proof of Corollary \ref{cornossb}

\medskip

\proof

Part (1).
Suppose for contradiction that the stated asymptotic self-similar blow-ups occur. Since the near maximal points are either at or within finite distance the side of $D$, we know that the limit of the blow-up solutions must live in $\mathbb{R}^2_+ \times (-\infty, 0]$. From Part (b) of the theorem, we know that $v_1=0$ and $v_2= c/(1-t)^\alpha$. But $v_2$ (corresponding to the vertical component of the velocity) is odd in the vertical variable by assumption. Hence $v_2=0$.
 Therefore $v \equiv 0$. This is a contradiction to the fact that $v$, as the blow-up limit $\tilde v$ with tilde dropped, is non-trivial.
 \medskip

Part (2).

In this part the spatial center of the blow-up point is the fixed boundary point $(1, 0)$ in the $r-x^{(3)}$ plane and we have ignored the angle $\theta$ due to axial symmetry. So we just need to replace the profile function $V$ by
\[
U=U(y) \equiv V(y+(1, 0)).
\]
The scaled solution $\td v_k$ is
\begin{equation}
\label{vktild3U}
\al
\tilde v_k &= \frac{Q_k^{-1}}{(\, 1-t_k-Q_k^{-1/\alpha} \,  \tilde t \, )^{\alpha} } V\left(\frac{Q_k^{-(1-\alpha)/\alpha} \tilde x}{[\, 1-t_k-Q_k^{-1/\alpha} \,  \tilde t \, ]^{1-\alpha} } +(1, 0) \right)\\
&=\frac{1}{(\, 1 - \tilde t \, )^{\alpha} } V\left(\frac{\tilde x}{(1-\tilde t)^{1-\alpha}} +(1, 0) \right)\\
&=\frac{1}{(\, 1 - \tilde t \, )^{\alpha} } U\left(\frac{\tilde x}{(1-\tilde t)^{1-\alpha}}  \right).
\eal
\end{equation} Notice that there is no more shift in the spatial direction as in part (1) where the center of blow-ups may move.
Using this and our assumption on $|\nabla V|$ and $|\nabla^2 V|$, we see,that
 $|y| |\nabla U|$, $|y|^2 |\nabla^2 U|$ are locally bounded functions with the reference point $y=0$.
 As in the proof of the theorem, after repeated differentiation of \eqref{vktild3U}, for any compact set $\Omega \subset D_k$ containing $\td x=0$ and $\td t \le 0$, we have
\be
\lab{vkttloc}
\Vert \partial_{\td t}  \tilde  v_k(\cdot, \tilde t) \Vert_{C^{0}(\Omega)}+ \Vert \partial_{\td t} \nabla \tilde  v_k(\cdot, \tilde t) \Vert_{C^{0}(\Omega)} +
 \Vert  \nabla \partial_{\td t} \td p_k(\cdot, \tilde t) \Vert_{C^{0}(\Omega)}+\Vert \partial_{\td t} \partial_{\td t} \tilde v_k(\cdot, \tilde t) \Vert_{C^{0}(\Omega)} \le C_6(\Omega),
\ee for a uniform constant $C_6(\Omega)$.
The arguments in part (b) of the theorem can be repeated verbatim.

Therefore, by Ascoli-Arzela theorem,  we still have local $C^{\gamma^-}_{x, t}$ convergence to a nontrivial limit solution $\td v$ in local $C^{2, 1, \gamma}_{x, t}$ space.
 Let $\omega$ be the (scalar) vorticity of $\tilde v$ again, which satisfies the transport equation \eqref{2dvtx}:
\[
\partial_{\tilde t} \omega + \tilde v \nabla \omega=0.
\]Since,  $\omega$ is bounded by assumption and vanishes at time $-\infty$  and $\td v$ has at most sub-linear growth, the flow lines are well defined and we can still conclude that $\omega =0$. Hence $\td v(\cdot, \td t)$ is harmonic. Now we drop  the tilde symbol and use the same notations as in the end of the proof of the Theorem, part (b). Since the center of the blow-ups is always at the boundary, in this case the limit domain $\tilde D$ is the half plane $\mathbb{R}^2_+$.
 Since $v_1=0$ at the boundary and $V=V(y)$ is sublinear and the center point of the blow-up is fixed, we know that $v_1$ is also sublinear at each time slice.   By the Liouville theorem again $v_1=0$. The divergence free condition $\partial_1 v_1 + \partial_2 v_2=0$ implies that $v_2=v_2(x_1, t)$.
Using $0= \omega = \partial_2 v_1 - \partial_1 v_2$, we see that $v_2=v_2(t)= c/(1-t)^\alpha$.
  But $v_2$ (corresponding to the vertical component of the velocity) is odd in the vertical variable by assumption. Hence $v_2=0$.
 Therefore $v=(v_1, v_2) \equiv 0$. This is a contradiction to the fact that $v$, as the blow-up limit $\tilde v$ with tilde dropped, is non-trivial to begin with.

The proof of the corollary is complete.

\qed

\medskip

\medskip

{\it Part II.}

Now we give a proof of Proposition \ref{prLH}.

\proof

(i). Suppose such ASSS exists. Since the error $W$ are local $C^{2, 2, \gamma}_{x, t}$  functions of lower order, without loss of generality we assume the error term is 0,  and also $T_0=1$.
 Following the proof of Corollary \ref{cornossb} part (2), we make a shift of the coordinates so the new origin $0$ is at $(1, 0)$ in the original $r-x^{(3)}$ coordinates. We see that the scaling and dimension reduction argument similar to \eqref{vktild3U} using the factor $Q_k=(1-t_k)^{-\alpha}$, which now converge to $0$, is still valid since
 \[
 x=Q_k^{-(1-\alpha)/\alpha} \tilde x + x_k, \quad t=
Q_k^{-1/\alpha} \tilde t + t_k,
 \]and therefore, due to $\alpha<0$, the scaled variables $\tilde x $ and $\tilde t$ are still blow ups of the original variables $x$ and $t$ respectively.  We have $\tilde v^\theta_k$ converges to $\tilde{v}^\theta$ in local $C^\gamma_x$ sense. However the components $\tilde{v}^r_k$ and $\tilde{v}^{(3)}_k$ converge to $0$ in local $C^\gamma_x$ sense, since they are $\beta$ order smaller than $v^\theta_k$. Indeed the k-th scaled solutions are
\be
\lab{LHsc1vtildk}
\al
&\tilde v^\theta_k=\tilde v^\theta_k(\tilde x, \tilde t) =\frac{1}{(1- \tilde t)^{\alpha}} \Theta \left( \frac{\tilde x}{(1- \tilde t)^{1-\alpha}} \right),  \\
&\tilde v^r_k=\tilde v^r_k(\tilde x, \tilde t) =\frac{Q^{-\beta/\alpha}_k}{(1- \tilde t)^{\alpha-\beta}} V^r \left( \frac{\tilde x}{(1- \tilde t)^{1-\alpha}} \right), \\
&\tilde v^{(3)}_k=\tilde v^{(3)}_k(\tilde x, \tilde t) =\frac{Q^{-\beta/\alpha}_k}{(1- \tilde t)^{\alpha-\beta}} V^{(3)} \left( \frac{\tilde x}{(1- \tilde t)^{1-\alpha}} \right).
\eal
\ee

The limiting solution is a one component solution of the Euler equation
\[
\tilde{v}=\tilde{v}^\theta e^\theta =(1-\tilde t)^{-\alpha} \Theta \left(\frac{\tilde x}{(1-\tilde t)^{1-\alpha}} \right) e_\theta, \quad \tilde t \le 0.
\]Notice that we do not have an independent $\tilde t$ variable in $\tilde v$ since the center of the blow up is fixed.
However, from the scaling invariant equation \eqref{eqgam} for $\Gamma = r v^\theta$, we deduce that the Euler equation \eqref{axi-EL} for $\tilde{v}$ collapses to
\[
\partial_{\tilde t} \tilde{v}=0.
\]Therefore in the unit ball, $\Theta=\Theta(y)$ is a homogeneous function of degree $\alpha/(1-\alpha)$, which must have a singularity at $0$. But $\Theta$ is locally bounded since we are scaling at the near maximal rate around the origin. This is a contradiction which proves part (i) of the proposition.

(ii). We follow the same set up as part (i) and perform the same scaling process.
Suppose the stated ASSS exists.
This time, however, the scaled solutions become
 \be
\lab{LHsc2vtildk}
\al
&\tilde v^\theta_k=\tilde v^\theta_k(\tilde x, \tilde t) =\frac{Q^{\beta/\alpha}_k}{(1- \tilde t)^{\alpha+\beta}} \Theta \left( \frac{\tilde x}{(1- \tilde t)^{1-\alpha}} \right),  \\
&\tilde v^r_k=\tilde v^r_k(\tilde x, \tilde t) =\frac{1}{(1- \tilde t)^\alpha} V^r \left( \frac{\tilde x}{(1- \tilde t)^{1-\alpha}} \right), \\
&\tilde v^{(3)}_k=\tilde v^{(3)}_k(\tilde x, \tilde t) =\frac{1}{(1- \tilde t)^\alpha} V^{(3)} \left( \frac{\tilde x}{(1- \tilde t)^{1-\alpha}} \right).
\eal
\ee Therefore $\tilde v^r_k$ and $\tilde v^{(3)}_k$ sub-converge in local $C^{2, 1, \gamma^-}_{x, t}$ sense to
 \[
 v^r_\infty =\frac{1}{(1- \tilde t)^\alpha} V^r \left( \frac{\tilde x}{(1- \tilde t)^{1-\alpha}}\right), \quad v^{(3)}_\infty =\frac{1}{(1- \tilde t)^\alpha} V^{(3)} \left( \frac{\tilde x}{(1- \tilde t)^{1-\alpha}}\right)
 \]on the half plane respectively; but $\tilde v^\theta_k$ diverge since $Q_k \to 0$ and $\beta>0$. However the normalized functions
\be
h_k \equiv Q^{-\beta/\alpha}_k \tilde v^\theta_k
\ee converge, in local $C^{2, 1, \gamma^-}_{x, t}$ sense to
\be
\lab{hinfty}
h_\infty \equiv \frac{1}{(1- \tilde t)^{\alpha+\beta}} \Theta \left( \frac{\tilde x}{(1- \tilde t)^{1-\alpha}} \right), \quad \tilde t \le 0.
\ee

We wish to take the limit for the equations satisfied by $v^r_k$, $v^{(3)}_k$ and $h_k$.
However, there is no uniform control of the pressure terms since $\tilde v^\theta_k$ become unbounded. Therefore, we will test the equations with compactly supported divergence free vector fields to eliminate the pressure. For this purpose, it is more convenient to use the following variables
$\tilde z =(\tilde z^{(1)}, \tilde z^{(2)})$
which are scaled version of $(r, x^{(3)})$, centered at $(1, 0)$.
\be
\lab{z1z2}
r=Q_k^{-(1-\alpha)/\alpha} \tilde z^{(1)} + 1,  \, x^{(3)}= Q_k^{-(1-\alpha)/\alpha} \tilde z^{(2)}, \quad t=
Q_k^{-1/\alpha} \tilde t + t_k.
\ee The volume element is the scaled axially symmetric one
\[
\tilde d \tilde z = r d\tilde z^{(1)} d\tilde z^{(2)}= (Q_k^{-(1-\alpha)/\alpha} \tilde z^{(1)} + 1)
d\tilde z^{(1)} d\tilde z^{(2)}.
\]The independent variables $\tilde x$ in  \eqref{LHsc2vtildk} are also replaced by $\tilde z$. To be precise, we should have given $\tilde z$ a subscript $k$. But for simplicity of presentation, we will skip this.

Similar to \eqref{eqtvro} and \eqref{eqtv3o} in the proof of the theorem,  the velocity functions in \eqref{LHsc2vtildk} satisfy the equations
\be
\lab{eqtvrz}
-( \tilde v^r_k \partial_{\tilde z^{(1)}} + \tilde
  v^{(3)}_k \partial_{\tilde z^{(2)}} ) \tilde v^r_k -\partial_{\tilde
  z^{(1)} }\tilde p_k - \partial_{\tilde t} \tilde v^r_k
+ \frac{( \tilde v^\theta_k)^2}{Q_k^{[(1/\alpha)-1]} r}=0.
\ee

\be
\lab{eqtv3z}
-( \tilde v^{r}_k \partial_{\tilde z^{(1)}} + \tilde
v^{(3)}_k \partial_{\tilde z^{(2)}} ) \tilde v^{(3)}_k
-\partial_{\tilde z^{(2)} }\tilde p_k -
\partial_{\tilde t} \tilde v^{(3)}_k=0.
\ee Notice that there is no error term since we are not converting to the rectangular coordinates before taking the limit.
 In the scaled variables, pick a smooth, compactly supported,  divergence free vector field $(\phi^r_k, \phi^{(3)}_k)$ which also vanishes on the boundary and use it as a test function in \eqref{eqtvrz} and \eqref{eqtv3z} respectively and take the sum. After integration in space, we see that all the pressure terms vanish:
\be
\lab{vth2kdk}
\al
&\int_{D_k} \frac{( \tilde v^\theta_k)^2}{Q_k^{[(1/\alpha)-1]} r} \phi^r_k \tilde d \tilde z= \int_{D_k} \left[ ( \tilde v^r_k \partial_{\tilde z^{(1)}} + \tilde
  v^{(3)}_k \partial_{\tilde z^{(2)}} ) \tilde v^r_k  + \partial_{\tilde t} \tilde v^r_k\right] \phi^r_k \tilde d \tilde z \\
&\qquad + \int_{D_k} \left[( \tilde v^{r}_k \partial_{\tilde z^{(1)}} + \tilde
v^{(3)}_k \partial_{\tilde z^{(2)}} ) \tilde v^{(3)}_k
+
\partial_{\tilde t} \tilde v^{(3)}_k \right] \phi^{(3)}_k \tilde d\tilde z.
\eal
\ee

 At this stage the domains are still finite, dilated cylinders and functions involved are still axially symmetric functions. We also choose the test vector fields $(\phi^r_k, \phi^{(3)}_k)$ in such a way that their supports are within a compact domain containing a fixed reference point and of fixed diameter in $D_k$ and their $C^{2, \gamma}_x$ norms are uniformly bounded.
 Multiplying \eqref{vth2kdk} by $Q^{-2(\beta/\alpha)+(1/\alpha)-1}_k$ on both sides and taking $k \to \infty$, using the convergence property of $h_k$, $\tilde v^r_k$ and $\tilde v^{(3)}_k$ we deduce, as long as $-2(\beta/\alpha)+(1/\alpha)-1>0$, that
\be
\lab{intTh2=0}
\int_{\mathbb{R}^2_+} \Theta^2 \left( \frac{(\tilde z^{(1)}, \tilde z^{(2)})}{(1- \tilde t)^{1-\alpha}} \right) \phi^r(\tilde z^{(1)}, \tilde z^{(2)}) d\tilde z^{(1)} d\tilde z^{(2)}=0, \qquad \tilde t \le 0.
\ee Let us explain why this is the case. Here we have used the dimension reduction idea in the proof of the theorem to infer that all the functions become 2 dimensional ones in the half plane, and axially symmetric divergence free property becomes the 2 dimensional divergence free property in the $\tilde z$ half plane.  Since the original radial variable $r=Q_k^{-(1-\alpha)/\alpha} \tilde z^{(1)} +1$ and $Q_k \to 0$ and $\alpha<0$, we see that $r \to 1$ locally. This also implies the convergence of volume elements to the 2 dimensional one, namely
\[
\tilde d \tilde z= r d\tilde z^{(1)} d\tilde z^{(2)} \to d\tilde z^{(1)} d\tilde z^{(2)}.
\]
In addition $(\phi^r, \phi^{(3)})$ is a local $C^{2, \gamma^-}_x$ limit of $(\phi^r_k, \phi^{(3)}_k)$ whose convergence can always be guarantied by suitable selection. Now that $(\phi^r, \phi^{(3)})$ is a smooth, compactly supported,  2 dimensional divergence vector field, we can write it as $\nabla^\bot f$ for some smooth function $f$. We can also choose the axially symmetric divergence free vector fields
\[
\left(\phi^r_k, \phi^{(3)}_k \right) = \left(-\frac{\partial_{\tilde z^{(2)}} f(\tilde z^{(1)}, \tilde z^{(2)} ) }{r}, \, \frac{\partial_{\tilde z^{(1)}} f(\tilde z^{(1)}, \tilde z^{(2)} ) }{r}\right), \quad r=Q_k^{-(1-\alpha)/\alpha} \tilde z^{(1)} +1,
\]so that the limit can reach any  $\nabla^\bot f$ with $f$ being smooth and compactly supported.
In particular $\phi^r= - \partial_{\tilde z^{(2)}} f$. Substituting this to \eqref{intTh2=0}, we see that $\Theta$ is independent of $\tilde z^{(2)}$. Since $\Theta$ is odd in the vertical variable $\tilde z^{(2)}$ by assumption, we conclude that $\Theta \equiv 0$. But $\Theta$ is non-trivial to begin with. We have reached a contradiction which proves the proposition if $-2(\beta/\alpha)+(1/\alpha)-1>0$.

(iii)  If $-2(\beta/\alpha)+(1/\alpha)-1=0$, we see that the test function argument in part (ii) imply that \eqref{eqtvrz}, \eqref{eqtv3z} and the scaled version of \eqref{eqgam} converge in local $C^\gamma_{x, t}$ sense to the Boussinesq equations for
$v_\infty = (v^r_\infty(z^{(1)}, z^{(2)}, t), v^{(3)}_\infty(z^{(1)}, z^{(2)}, t))$ and $h_\infty=h_\infty(z^{(1)}, z^{(2)}, t)$ in the left half plane $z=(z^{(1)}, z^{(2)}) \in \mathbb{R}^2, \, z^{(1)} \le 0$. i.e.
\be
\lab{bouseq}
\begin{cases}
\al
&v_\infty \nabla v_\infty + \nabla p + \partial_t v_\infty - (h^2_\infty, 0)=0, \quad t \le 0,\\
&v_\infty \nabla h_\infty + \partial_t h_\infty=0,  \quad div \, v_\infty=0, \quad v^r_\infty(0, z^{(2)})=0.
\eal
\end{cases}
\ee Here we are using the same notations for gradient, pressure and the time variable, and drop the tilde in the $z$ variables  in order to save symbols. This completes the proof of the proposition. \qed
\medskip

Next we state and prove a non-existence result for the case $-2(\beta/\alpha)+(1/\alpha)-1=0$, i.e.  $2 \beta = 1- \alpha$. See the remark after the proof the reason of having  many conditions in the proposition.

\begin{proposition}
\lab{prLH00}
Consider the Boussinesq equations in the upper half of the $z=(z^{(1)}, z^{(2)})$ plane:
\be
\lab{bouseq02}
\begin{cases}
\al
&v \nabla v + \nabla p + \partial_t v - (0, h^2)=0, \quad z^{(2)} \ge 0, \, t \le 0,\\
&v \nabla h+ \partial_t h=0,  \quad div \, v=0.
\eal
\end{cases}
\ee
 Suppose $-2(\beta/\alpha)+(1/\alpha)-1=0$ with $\alpha<0$ and $0< \beta <|\alpha|$. There does not exist self-similar blow-up solutions
 \[
 \al
 v&=v(z, t)=(v^{(1)}(z, t), v^{(2)}(z, t))\\
 &=\frac{1}{(1-  t)^\alpha} \left( V^{(1)} \left( \frac{z}{(1-  t)^{1-\alpha}} \right), \, V^{(2)} \left( \frac{z}{(1-  t)^{1-\alpha}} \right) \right),\\
 h&=h(z, t)=\frac{1}{(1-  t)^{\alpha+\beta}} H \left( \frac{z}{(1-  t)^{1-\alpha}} \right)
 \eal
 \]to \eqref{bouseq02} with the properties given below.

  (a). The profiles are nontrivial  local $C^{2, \gamma}_x$  functions whose gradients and Hessian are uniformly bounded; $V^{(1)}$ and $H$ are odd in $z^{(1)}$, $V^{(2)}$ is even in $z^{(1)}$; $V^{(2)}=0$ at the base of the half plane. $H^2(z)>0$ when $|z^{(1)}| \to \infty$.

  (b). $|V^{1}(z)|+|V^{2}(z)| \le C (1+|z|)^{\alpha/(\alpha-1)}$ and $|H(z)| \le C (1+|z|)^{(\alpha+\beta)/(\alpha-1)}$ for a positive constant $C$;

  (c).  When $|z|$ is sufficiently large $V^{(2)}(z) \ge 0$.

   (d). Let $W=W(z)$ be the profile of $\omega=\partial_{z^{(2)}} v^{(1)} -\partial_{z^{(1)}} v^{(2)}$.  There is a large constant $l_0>0$, When $z^{(1)}, z^{(2)} \in [0, l_0]$, there is one $C^1$ curves $ z^{(1)}=L(z^{(2)})$ along which $\partial_{z^{(1)}} W(z)=0$. Also $\partial_{z^{(1)}} W(z) < 0$ when $z^{(1)} > L(z^{(2)})$, $\partial_{z^{(1)}} W(z) >0$ when $z^{(1)} \le L(z^{(2)})$. .

   (e).
    $\partial_{z^{(1)}} H^2(z) \ge 0$ near the curve in (d). Also $|\partial_{z^{(1)}} H^2/\partial_{z^{(1)}} W|$ is sub-linear in $z$ when $|z|$ is large.

\end{proposition}
 \proof
Notice that we have made a rotation and  reflection  so that the problem in the left half of the plane becomes that in the upper half.  The former is derived from suitable blow up from the ASE problem at the vertical cylindrical boundary.  The latter is in line with relevant literature.
 Let
\be
\lab{omprofil2}
 \omega=\omega(z, t)=\partial_{z^{(2)}} v^{(1)}- \partial_{z^{(1)}} v^{(2)}= \frac{1}{(1-  t)} W \left( \frac{z}{(1-  t)^{1-\alpha}} \right)
\ee be the scalar vorticity. Note the sign convention for $\omega$ is the opposite to that in the standard text book. It is the same as in \cite{CH1} but opposite to that of \cite{WLG-SB}.

 Then $\omega$ satisfies the equation
 \be
\lab{2dbwtx2}
\partial_{ t} \omega +  v \nabla \omega + \partial_{z^{(1)}} h^2=0.
\ee  Consider $\omega $ in the infinite strips in the $z$ plane:
\[
S_t \equiv \{ (z^{(1)}, z^{(2)} ) \, | \, z^{(1)} \ge 0,   \, 0 \le z^{(2)} \le l_0 (1-t)^{1-\alpha} \}.
\]We rewrite equation \eqref{2dbwtx2} as
\be
\lab{2dbwtx22}
\partial_{ t} \omega +
\left[ v^{(1)} + \left(\partial_{z^{(1)}} h^2/\partial_{z^{(1)}} \omega \right) \right]\partial_{z^{(1)}} \omega  +  v^{(2)} \partial_{z^{(2)}} \omega =0.
\ee Notice that the term $\partial_{z^{(1)}} h^2/\partial_{z^{(1)}} \omega$ may have  singularities when $\partial_{z^{(1)}} \omega=0$. But under assumption (d), at each time level $t \le 0$, the singularities can happen only along a $C^1$ curve:
\[
z^{(1)}/(1-t)^{1-\alpha}=L(z^{(2)}/(1-t)^{1-\alpha}),  \quad z^{(2)} \in [0, l_0 (1-t)^{1-\alpha}],
\] which is called $l$.

Starting from a point $z \in S_0$ and $t=0$, we claim that the backward flow line $c=c(t)=(z^{(1)}(t), z^{(2)}(t)),  t \le 0$ of the vector field
\be
\lab{xlcc}
\left(v^{(1)} + (\partial_{z^{(1)}} h^2/\partial_{z^{(1)}} \omega ), \, v^{(2)} \right)
\ee can either be extended continuously to $t=-\infty$ or it will hit the left side $z^{(1)}=0$.

Taking this claim for granted, we will show that the conclusion of the proposition is true.
First we show that $c(t) \in S_t$ unless $c(t_1)$  hits the left side at some $t_1>t$. The reason is that $v^{(2)}=0$ at the base $z^{(2)}=0$, so that $c(t)$ can not cross the base. If $c(t)$ hits  the top, then $z^{(2)}(t)=l_0 (1-t)^{1-\alpha}$, assumption (c) infers that $v^{(2)} =(1-  t)^{-\alpha}V^{(2)}(z^{(1)}(t)/(1-t)^{1-\alpha},  l_0) \ge 0$. Therefore the flow line can not cross the top of $S_t$ backward in time. Hence the only way for the backward flow line to escape $S_t$ is through the left side where $\omega=0$ by the oddness assumption. In this case $\omega=0$ along the whole flow line.
Consequently, for all $ t \le 0$
\[
|\omega(z, 0)| \le \sup_{S_t} |\omega(\cdot, t)| \to 0, \quad t \to -\infty.
\]Here we have used \eqref{omprofil2}. Hence $\omega(z, 0)=0$, for all $z \in S_0$. From \eqref{2dbwtx2}, this infers
$h^2$ is independent of $z^{(1)}$. But $h=0$ at the left side, which shows that $h^2=0$ on $S_0$. This is a contradiction with part of the assumption (a) that $h^2(z, 0)=H^2(z)$ is positive as $|z| \to \infty$.  So the proposition is true if the claim is true.

Now we give a proof of the claim. The main point is that the flow line can not cross the $C^1$ curve $l$ as long as the assumptions on $\omega, h^2$ hold.
Observe that the flow lines of \eqref{xlcc} satisfy, except at the singularities, the equations
\[
\frac{d z^{(1)}(t)}{d t} = [v^{(1)}+ (\partial_{z^{(1)}} h^2/\partial_{z^{(1)}} \omega )],  \quad
\frac{d z^{(2)}(t)}{d t} = v^{(2)}(z^{(1)}, z^{(2)}, t).
\] By the assumption on $v$, $h^2$ and $w$, the pertinent vector fields are sub-linear in space, so the flow lines will continue backward in time until it hits a singularity on $l$.

Let us make the change of variables
\[
y_1=z^{(1)} - L(z^{(2)}), \quad y_2 = z^{(2)}.
\]Then the equations for $y_1$ and $y_2$ become
\[
\al
&\frac{dy_1(t)}{d t} = [v^{(1)}+ (\partial_{z^{(1)}} h^2/\partial_{y_1} \omega )]
 - L'(y_2) v^{(2)}(y_1+L(y_2), y_2, t),\\
&\frac{dy_2(t)}{d t} = v^{(2)}(y_1+L(y_2), y_2, t).
\eal
\]Also the curve $l$ becomes a segment of $y_1=0$.

Now pick a point to the right of the singularity, we see that the flow line ending at this point can not cross the singularity backward. The reason is that $v^{(1)}- (\partial_{z^{(1)}} h^2/\partial_{y_1} \omega )$ approaches $-\infty$. This can be seen if one compute directly that for the
variable
\be
\lab{tldz}
\tilde{z}=(\tilde{z}^{(1)}, \tilde{z}^{(2)}) =(z^{(1)}/(1-t)^{1-\alpha}, z^{(2)}/(1-t)^{1-\alpha}),
\ee we have
\be
\lab{eqtldz}
(1-t) \frac{d\tilde{z}^{(1)}}{dt}= \left[V^{(1)}(\tilde{z})
+ \frac{\partial_{\tilde{z}^{(1)}} H^2}{\partial_{\tilde{z}^{(1)}} W}(\tilde{z}) + (1-\alpha)\tilde{z}^{(1)}  \right].
\ee By assumption (d) and (e), it is clear that $\frac{d\tilde{z}^{(1)}}{dt}$ approaches $-\infty$ if $\tilde{z}$ approaches the singularity (zero of $\partial_{\tilde{z}^{(1)}} W)$ from the right and it approaches $\infty$ if $\tilde{z}$ approaches the singularity from the left.
Therefore, being a small perturbation,
\[
\frac{d(y^{(1)}/(1-t)^{1-\alpha})}{dt} \le 0
\]in a small right neighborhood of the singularity. So $y^{(1)}/(1-t)^{1-\alpha}$ will move away from the singularity along the backward flow line. Similarly, if the end point of the flow line is to the left of the singularity, the flow line will also not cross the singularity backward.

 Converting back to the original $z$ variable, we know that its flow lines can be extended continuously too. This conversion is always possible due to the inverse function theorem since the Jacobian between $(y_1, y_2)$ and $(z^{(1)}, z^{(2)})$ is $1$.   We hereby finish the proof of the claim, except at the singularities (curves given in (d)). This implies $\omega=0$ except at the curves. By continuity, $\omega=0$ everywhere.   The proposition is proven.
\qed

{\remark (i).  This seems to include the expected self similar solution (SSS) of the Boussinesq equation in \cite{WLG-SB} and the Euler equation outside a cylinder. See p2 there for the explanation on converting the Euler equation to \eqref{bouseq02}. Note the vorticity there is the negative of the vorticity here. Also it seems that the figure for the function $\Phi$
on p2 there should be reflected across the horizontal axis since it is the horizontal i.e. $y_1$ derivative of a square function. What is computed in \cite{WLG-SB} appears to be the same interior blow up scenario as in \cite{CH1}.

 Actually one can just work on the base $z^{(2)}=0$ to quickly reach the conclusion $\omega =0$ and $h^2=0$ on the base. The reason is $v^{(2)}=0$ on the base so  $z^{(2)}$ stays $0$ on the flow line. So \eqref{eqtldz} reduces to
\be
\lab{eqtldz0}
(1-t) \frac{d\tilde{z}^{(1)}}{dt}= \left[V^{(1)}(\tilde{z}^{(1)}, 0)
+ \frac{\partial_{\tilde{z}^{(1)}} H^2}{\partial_{\tilde{z}^{(1)}} W}(\tilde{z}^{(1)}, 0) + (1-\alpha)\tilde{z}^{(1)}  \right].
\ee

All the assumptions of the proposition can be made only on the base and there is no need to change variables. This already induces a contradiction with assumption (a). What we did right before the remark is to reach a stronger conclusion that $\omega=0$ and $h^2=0$ everywhere. This is true since $W$ is $0$ on the vertical axis and $0$ at infinity. So, for any $z^{(2)}>0$,  there is always a $z^{(1)}>0$ such that $\partial_{z^{(1)}} W=0$ at
$(z^{(1)}, z^{(2)})$. So we can extend the proof of the proposition to the whole upper plane.

(ii). Notice the sign in front of $\partial_{z^{(1)}} h^2$ in \eqref{2dbwtx2} is positive.
In \cite{CH1} the sign is negative. This would not matter in general by changing $h^2$ to $-h^2$. However as the Boussinesq equation arrives as a scaling limit from ASE, the difference in sign becomes important since $-h^2$ decreases to negative infinity as $z^{(1)} \to \infty$.

(iii). Some relaxation of the condition on $W$ is possible. For example one can assume that on the base $z^{(2)}=0$, $W$ can have finitely many local maximum and minimums as long as the last one on the right is a local maximum.  Then flow lines ending at the right of the maximal point can be extended to $t=-\infty$. This shows $\omega=0$ in a line segment which is also a contradiction.}

\medskip

Next we state and prove a non-existence result for the case $-2(\beta/\alpha)+(1/\alpha)-1=0$, i.e.  $2 \beta = 1- \alpha$  when the velocity/vorticty equation has a different sign in front of $h^2$. This case corresponds to potential vorticity blow up of ASE at the inside boundary of a cylinder. There appear to be too many conditions at the first glance. But they are tailored for some expected blow up scenario.
The main point is that under this scenario, the profile of the vorticity $W$ can not have a saddle point or local maximum point outside a compact set.

\begin{proposition} (no interior saddle, local maximum.)
\lab{prLH0}
Consider the Boussinesq equations in the upper half of the $z=(z^{(1)}, z^{(2)})$ plane:
\be
\lab{bouseq0}
\begin{cases}
\al
&v \nabla v + \nabla p + \partial_t v + (0, h^2)=0, \quad z^{(2)} \ge 0, \, t \le 0,\\
&v \nabla h+ \partial_t h=0,  \quad div \, v=0.
\eal
\end{cases}
\ee
 Suppose $-2(\beta/\alpha)+(1/\alpha)-1=0$ with $\alpha<0$ and $0< \beta <|\alpha|$.  Suppose the functions
 \be
 \lab{asssvh}
 \al
 v&=v(z, t)=(v^{(1)}(z, t), v^{(2)}(z, t))\\
 &=\frac{1}{(1-  t)^\alpha} \left( V^{(1)} \left( \frac{z}{(1-  t)^{1-\alpha}} \right), \, V^{(2)} \left( \frac{z}{(1-  t)^{1-\alpha}} \right) \right),\\
 h&=h(z, t)=\frac{1}{(1-  t)^{\alpha+\beta}} H \left( \frac{z}{(1-  t)^{1-\alpha}} \right)
 \eal
 \ee are self-similar solutions to \eqref{bouseq0} with the following properties.
The profiles are local $C^{2, \gamma}_x$  functions whose gradients and Hessian are uniformly bounded; $V^{(1)}$ and $H$ are odd in $z^{(1)}$, $V^{(2)}$ is even in $z^{(1)}$; $V^{(2)}=0$ at the base of the half plane.

(i). Then no such SSS $v$ or ASSS with $v$ as scaling limit exists under the further conditions  (a)-(d).

(a).  In the open first quadrant, $V^{(1)} < 0$, $V^{(2)} > 0$, and $\partial_{z^{(1)}} H^2 >0$ except on the vertical axis, and $V$ is sub-linear near infinity.
Let $W=W(z)$ be the profile of $\omega=\partial_{z^{(2)}} v^{(1)} -\partial_{z^{(1)}} v^{(2)}$,  $\sup W= \sup |W| >0$ in the first quadrant.

(b). Near infinity, for $r=|z|, \, \theta = \tan (z^{(2)}/z^{(1)})$,
\[
W(z)=W(r, \theta)= \frac{f(r, \theta)}{r^{1/(1-\alpha)}}, \quad V_r(r, \theta)= \frac{g(r, \theta)}{r^{\alpha/(1-\alpha)}},
\]where $f(r, \theta)$ and $g(r, \theta)$ are bounded functions.

 (c). There is an angle $\theta_1>0$ and $\theta_2 \in (\theta_1, \pi/2)$ such that
\[
\sup_{r>0} W(r, \theta_2) > \sup_{r>0} W(r, \theta_1).
\]

(d). $W \ge 0$ and
$W$ reaches strict maximum value $W(z_0)>0$ in the sector $D_{\infty, \theta_1, \theta_2}=\{(r, \theta)  \, | \,  r \ge 0, \, \theta \in [\theta_1, \theta_2]\}$.
i.e. for any small $\epsilon>0$,
$W(z_0)>\sup_{B^c(z_0, \epsilon) \cap D_{\infty, \theta_1, \theta_2}} W$.
Moreover,
\[
W(z_0)- \partial_{z^{(1)}} H^2(z_0)>0.
\]

(ii). The same conclusion still holds under condition (a) and the next condition:

(b').  Let $P_i=(a_i, b_i), i=1, 2, 3, 4$ be four points in the first quadrant which form a rectangle $D$ whose sides are parallel to the coordinate axis. i.e.  $a_1=a_4,  b_1=b_2, a_2=a_3,  b_3=b_4$, $a_2>a_1$ and $b_4>b_1$ so that $P_1$ is at the lower left corner. $W \ge 0$ and a strict maximum of $W$ in $D$ occurs either in the interior of $D$ or the interior of the union of the upper side $P_3P_4$ and the right side $P_2P_3$.
Also $
W \ge \partial_{z^{(1)}} H^2
$ in $D$ and $V^{(1)} + (1-\alpha) z^{(1)}>0$ except when $z^{(1)}=0$.

\end{proposition}
 \proof

Part (i).  Suppose conditions (a), (b), (c) and (d) hold and such SSS $v$ exists.

One equation for the profiles in the upper $z$ plane reads
\be
\lab{eqprf}
V \nabla W + (1- \alpha) z \nabla W + W - \partial_{z^{(1)}} H^2=0.
\ee Taking $l>1$ and the two angles $\theta_1, \theta_2 \in (0, \pi/2)$ and pick the sectorial domain
\[
D_{l, \theta_1, \theta_2}=B^I(0, l) \cap \{\theta_1<\theta<\theta_2 \}
\]where $B^I(0, l)$ is the disk with radius $l$ in the first quadrant. Using $W^p$ with a large $p>1-2 \alpha$ as a test function on $D_{l, \theta_1, \theta_2}$,  we deduce, using the divergence free property of $V$ and boundary conditions, that
\be
\lab{w123}
\al
&0=\int_{D_{l, \theta_1, \theta_2}} \left[  \left( 1- \frac{2(1-\alpha)}{p+1} \right)W- \partial_{z^{(1)}} H^2 \right] \, W^p dz\\
&\qquad + \frac{1}{p+1} \int^{\theta_2}_{\theta_1} V_r l \,  W^{p+1} (l, \theta)  d\theta
+ \frac{1-\alpha}{p+1} \int^{\theta_2}_{\theta_1} l^2 W^{p+1} (l, \theta) d\theta\\
&\qquad+
\frac{1}{p+1} \int^{l}_0 V \cdot (-\sin \theta_2, \cos \theta_2) \,  W^{p+1} (r, \theta_2)  dr \\
&\qquad-\frac{1}{p+1} \int^{l}_0 V \cdot (-\sin \theta_1, \cos \theta_1) \,  W^{p+1} (r, \theta_1)  dr.
\eal
\ee Here $V_r = V \cdot (\cos \theta, \sin \theta)$.

Since $p>1-2 \alpha$, using the assumption on the behavior of $V$ and $W$ at infinity in (b), we find, after taking $l \to \infty$, that
\be
\lab{winfty}
\al
&0=\int_{D_{\infty, \theta_1, \theta_2}} \left[  \left( 1- \frac{2(1-\alpha)}{p+1} \right)W- \partial_{z^{(1)}} H^2 \right] \, W^p dz\\
&\qquad+
\frac{1}{p+1} \int^{\infty}_0 V \cdot (-\sin \theta_2, \cos \theta_2) \,  W^{p+1} (r, \theta_2)  dr\\
&\qquad -\frac{1}{p+1} \int^{\infty}_0 V \cdot (-\sin \theta_1, \cos \theta_1) \,  W^{p+1} (r, \theta_1)  dr.
\eal
\ee By our assumption $V^{(1)} < 0$ and $V^{(2)} > 0$ in  the open first quadrant, the function $V \cdot (-\sin \theta, \cos \theta)= -V^{(1)} \sin \theta + V^{(2)} \cos \theta$ is positive in the open first quadrant. Therefore,
\[
\lim_{ p \to \infty}  \left(\int^{\infty}_0 V \cdot (-\sin \theta, \cos \theta) \,  W^{p+1} (r, \theta)  dr \right)^{1/(p+1)} = \sup_{r>0} W(r, \theta).
\]Here we also used the condition that $V$ is sub-linear  and $W$ decays near infinity. By our assumption
\[
\sup_{r>0} W(r, \theta_2) > \sup_{r>0} W(r, \theta_1),
\]for $p$ sufficiently large,  we have
\be
\lab{vth12}
\int^{\infty}_0 V \cdot (-\sin \theta_2, \cos \theta_2) \,  W^{p+1} (r, \theta_2)  dr-
 \int^{\infty}_0 V \cdot (-\sin \theta_1, \cos \theta_1) \,  W^{p+1} (r, \theta_1)  dr >0.
\ee

By assumption (c), there are fixed small numbers $\delta$ and $\epsilon$ such that
\[
 \left( 1- \delta \right)W- \partial_{z^{(1)}} H^2 >0
\]in the region  $B(z_0, \epsilon) \cap D_{\infty, \theta_1, \theta 2}$.
This leads to, for large $p$, that
\[
 \left( 1- \frac{2(1-\alpha)}{p+1} \right)W- \partial_{z^{(1)}} H^2>0
\]in $B(z_0, \epsilon) \cap D_{\infty, \theta_1, \theta 2}$.
Hence, we deduce
\[
\al
&\lim_{p \to \infty} \left( \int_{B(z_0, \epsilon) \cap D_{\infty, \theta_1, \theta 2}} \left[  \left( 1- \frac{2(1-\alpha)}{p+1} \right)W- \partial_{z^{(1)}} H^2 \right] \, W^p dz \right)^{1/(p+1)}\\
&= \sup_{B(z_0, \epsilon) \cap D_{\infty, \theta_1, \theta 2}} W=W(z_0).
\eal
\]In the region $B^c(z_0, \epsilon) \cap D_{\infty, \theta_1, \theta 2}$, we have,
\[
\al
&\lim_{p \to \infty} \left( \int_{B^c(z_0, \epsilon) \cap D_{\infty, \theta_1, \theta 2}} \left|  \left( 1- \frac{2(1-\alpha)}{p+1} \right)W- \partial_{z^{(1)}} H^2 \right| \, W^p dz \right)^{1/(p+1)} \\
&\le \sup_{B^c(z_0, \epsilon) \cap D_{\infty, \theta_1, \theta 2}} W<W(z_0)
\eal
\]by assumption (d). Note the zeroes of $\left( 1- \frac{2(1-\alpha)}{p+1} \right)W- \partial_{z^{(1)}} H^2$ may fall on the maximum points of $W$.
Hence for $p$ sufficiently large, we see from the last two lines that
\be
\lab{Dwp+1>0}
\int_{D_{\infty, \theta_1, \theta_2}} \left[  \left( 1- \frac{2(1-\alpha)}{p+1} \right)W- \partial_{z^{(1)}} H^2 \right] \, W^p dz>0
\ee

Thus  \eqref{winfty} is a contradiction by \eqref{vth12} and \eqref{Dwp+1>0}.
We hereby finish the proof of the proposition under  conditions (a), (b), (c) and (d).

Part (ii).  Suppose conditions (a), (b') hold and the conclusion is false. i.e. such SSS exists.

Let $z=z(s)=( z^{(1)}(s), z^{(2)}(s)) $ be an integral curve (flow line) of the vector field $V + (1-\alpha) z$, i.e. \be
\lab{zflowline}
\frac{d z(s)}{d s} = V(z(s)) + (1-\alpha) z(s)=(V^{(1)} + (1-\alpha) z^{(1)}, V^{(2)} + (1-\alpha) z^{(2)}).
\ee This and
equation \eqref{eqprf}  tell us
\be
\lab{wrect}
\frac{d W(z(s))}{d s} + (W-\partial_{z^{(1)}} H^2)(z(s))=0
\ee From \eqref{zflowline} and the conditions $V^{(1)} + (1-\alpha) z^{(1)}>0$ and $V^{(2)}>0$ in $D$, we see that  $\frac{d z(s)}{d s}$ which is the tangent vectors of the flow line forms an acute angle with respect to the positive $z^{(1)}$ axis.

Suppose the strict maximum value  of $W$, say $M>0$, occurs in the interior point of $D$, say $z_0$. Then the backward flow line $z(s)$ ending at $z_0$, i.e. $z(0)=z_0$ will stay in $D$ for a short time $s<0$. By \eqref{wrect},  $W(z(s)) \ge W(z_0)$. This is a contradiction with $W(z_0)$ being a strict maximum value.

Next, suppose the strict maximum value  of $W$, say $M>0$, occurs in the interior point $z_0$ of the union of the upper side $P_3P_4$ and the right side $P_2P_3$.
Then, due to acuteness of the angle of $\frac{d z(s)}{d s}$,  the backward flow line $z(s)$ ending at $z_0$  will still stay in $D$ for a short time $s<0$. By \eqref{wrect},  $W(z(s)) \ge W(z_0)$. This is a contradiction with $W(z_0)$ being a strict maximum value.
This completes the proof of the proposition.
\qed
\medskip

{\remark

  The expected SSS of the Boussinesq equation in \cite{CH1} (c.f. Sec. 5 and Figure 1, p9)  appear  to satisfy all the assumptions in the proposition if that figure is accurate. See also Remark \ref{relocalsadd} below. Note that condition (c) is satisfied due to the dip in the ridge of $W$ in Figure 1, and our sector cuts away the neighborhood where
    $W-\partial_{z^{(1)}} H^2=0$
     except at the vertical axis. Condition (d) is a local condition since the growth of $H^2(z)$ is of order $|z|^{2(\alpha+\beta)/(\alpha-1)}$ and the decay of $\partial_{z^{(1)}} H^2(z)$ is of the order $|z|^{-2/(1-\alpha)}$, which is faster than that of $W$. Inside a compact domain, we work under the assumption and belief that the figure and data in Figure 1 are accurate.  The angle  $\theta_1 $ is chosen so that the ray $\theta=\theta_1$ goes under the saddle  (the lowest point of the ridge of $W$) and $\theta_2$ is a little larger. There is a rise of W along the yellow ridge  in  Figure 1 after the lowest point of the saddle, in the direction of positive y  axis ($z^{(2)}$ axis here).
That is why $sup_r W(r, \theta_2)$ is larger than $sup_r W(r, \theta_1)$. The strict maximum value of $W$ in the sector $D_{\infty, \theta_1, \theta_2}$ is a little larger than $0.5$, which occurs on the ray $\theta=\theta_2$. At this point, $\partial_{z^{(1)}} H^2$ is around $0.3$. Therefore condition (d) is satisfied.

   The order of growth for the velocity $V$ is around $2/3$ and the decay for the vorticity $W$ is around $-1/3$ at $\infty$ and the order of growth of $H^2$ is around $1/3$ and $\partial_{x^{1}} H^2$ has a decay around $-2/3$. c.f. Sec.5.1 and Sec.6.4.2 \cite{CH1}. Note that it was not explicitly stated in \cite{CH1} that exact self similar solution (SSS) exists but certain approximate SSS exists which would expectedly converges to an exact SSS.
}

{\remark
\lab{relocalsadd}
There is a suggestion that the saddle shape in Figure 1 of \cite{CH1} is a result of perspective in the picture and is not real. But the integral argument in the proposition can be localized as follows. Consider the sectorial domain
\[
D=D_{l_1, l_2, \theta_1, \theta_2}=\{ (r, \theta) \, |  l_1<r<l_2, \quad \theta_1<\theta<\theta_2 \}.
\] Using $W^p$ with a large $p>1-2 \alpha$ as a test function on $D_{l, \theta_1, \theta_2}$,  we deduce, using the divergence free property of $V$ and boundary conditions, that
\be
\lab{w123loc}
\al
&0=\underbrace{\int_{D} \left[  \left( 1- \frac{2(1-\alpha)}{p+1} \right)W- \partial_{z^{(1)}} H^2 \right] \, W^p dz}_{T_1}\\
&+ \underbrace{ \frac{1}{p+1} \int^{\theta_2}_{\theta_1} V_r l_2 \,  W^{p+1} (l_2, \theta)  d\theta
+ \frac{1-\alpha}{p+1} \int^{\theta_2}_{\theta_1} l^2_2 W^{p+1} (l_2, \theta) d\theta}_{T_2}\\
& - \big[ \underbrace{ \frac{1}{p+1} \int^{\theta_2}_{\theta_1} V_r l_1 \,  W^{p+1} (l_1, \theta)  d\theta
+\frac{1-\alpha}{p+1} \int^{\theta_2}_{\theta_1} l^2_1 W^{p+1} (l_1, \theta) d\theta}_{T_3} \big]+ \\
&
\underbrace{\frac{1}{p+1} \int^{l_2}_{l_1} V \cdot (-\sin \theta_2, \cos \theta_2) \,  W^{p+1} (r, \theta_2)  dr }_{T_4}
- \underbrace{\frac{1}{p+1} \int^{l_2}_{l_1} V \cdot (-\sin \theta_1, \cos \theta_1) \,  W^{p+1} (r, \theta_1)  dr}_{T_5}.
\eal
\ee If the strict maximum of $W$ in $D$ occurs at the interior of the union of the  upper edge $\theta=\theta_2$ and the upper arc $r=l_2$, then we can rearrange the above as
\[
T_1+T_2+T_4=T_3+T_5.
\]By the same argument as in the proposition, as $p \to \infty$, the left side dominates the right side and  we still reach a contradiction.

From the data file "solu.w\{1, 1\}" in reference [1] of \cite{CH2} in the folder file
\[
"Steady\_state\_pertb720\_Nlevcor4.mat"
\] and sub folder "solu",
there are several rectangles where the above local condition in the sectors or condition (b') seem to hold if one uses piecewise affine interpolation between mesh points. For instance, the rectangle given by the $x-y$ ($z^{(1)}-z^{(2)}$) mesh points $(14, 719)$,  $(15, 719)$, $(14, 718)$ and $(15, 718)$. The strict maximum value of $W=9.02 \times 10^{-17}$ occurs at the upper right corner. Another example is the rectangle given by the mesh points $(19, 710),  (19, 712), (21, 710), (21, 712)$. Here interior maximum value $2.3461 \times 10^{-16}$ is reached at $(20, 711)$. We comment that in the far field, at small scales, rectangles are close to the sectorial domains.

In addition, at mesh point $(36, 720)$, we see from the $solu.w\{1, 1\}$ file that
\be
\lab{w-17}
W=8.843970087044398 \times 10^{-17}.
\ee From the files $solu.u1\{1, 2\}$ and $solu.u2\{2, 1\}$, we find
\[
\partial_{z^{(2)}} V^{(1)}= -2.850898611910781 \times 10^{-19}, \qquad \partial_{z^{(1)}} V^{(2)}=1.023527319550463 \times 10^{-5}
\]respectively. Therefore
\be
\lab{w-5}
W=\partial_{z^{(2)}} V^{(1)} - \partial_{z^{(1)}} V^{(2)} \approx -1.023 \times 10^{-5}.
\ee  From \eqref{w-17} and \eqref{w-5} we see the difference in the order $10^{-5}$ with opposite sign. This is beyond the accuracy of $10^{-7}$ as stated on p64 of \cite{CH1} for the approximate steady state $W$. Notice that the weighted pointwise norm on p8 and p137 there will only increase the error near this mesh point since $p_{6, 5} |z|^{1/7}$ on p137, which is part of the weights, is quite large.

In addition, a negative value for $W$ in the open first quadrant will lead to a contradiction as explained in the next remark.
}

{\remark Proposition \ref{prLH0} does not rule out all ASSS.
In general, suppose, when $z^{(1)}>0$, both $V^{(1)}(z^{(1)}, 0) + (1-\alpha) z^{(1)} >0$ and
  $\partial_{z^{(1)}}  H^2(z^{(1)}, 0) \ge 0$. Let $W=W(z)$ be the profile of $\omega$. If there is a point $z^{(1)}_0>0$ such that $W(z^{(1)}_0, 0)=0$, then $\omega =0$, $h^2=0$ at the base $z^{(2)}=0$. In particular, we have the following observation: if $W$ has one extra $0$ beside the origin, there does not exist a self similar solution such that $h^2$ is positive in the positive horizontal axis.

Here goes the proof. Let
\be
\lab{omprofil}
 \omega=\omega(z, t)=  \partial_{z^{(2)}} v^{(1)} -\partial_{z^{(1)}} v^{(2)}= \frac{1}{(1-  t)} W \left( \frac{z}{(1-  t)^{1-\alpha}} \right)
\ee be the scalar vorticity. Note the sign convention for $\omega$ is the opposite to that in the standard text book. One reason is that it is the suitable limit of the angular vortivity in the axially symmetric case, where the flow is in the finite cylinder. There $\omega_\theta= \partial_{x^{(3)}} v^r - \partial_r v^{(3)}$. After the above rotation and reflection,  $x^{(3)}$ axis become $z^{(1)}$ axis and $r$ axis becomes $z^{(2)}$ axis.  Then $\omega$ satisfies the equation
 \be
\lab{2dbwtx}
\partial_{ t} \omega +  v \nabla \omega - \partial_{z^{(1)}} h^2=0,
\ee $h^2$ satisfies
 \be
\lab{2dbh2tx}
\partial_{ t} h^2 +  v \nabla  h^2=0.
\ee Using \eqref{asssvh} and \eqref{omprofil}, we see that the following equations hold at the base $z^{(2)}=0$.
\be
\lab{eqWH2}
\al
&\left(V^{(1)}  + (1-\alpha) z^{(1)} \right) \partial_{z^{(1)}} W + W - \partial_{z^{(1)}} H^2=0,\\
&\left(V^{(1)}  + (1-\alpha) z^{(1)} \right) \partial_{z^{(1)}} H^2 + (1+\alpha) H^2=0.
\eal
\ee Under our assumption, this system has explicit solutions at the base $z^{(2)}=0$:
\[
\al
&H^2(z^{(1)}, 0)= C \exp (-(1+\alpha) \int (V^{(1)}  + (1-\alpha) z^{(1)} )^{-1} dz^{(1)}),\\
&W(z^{(1)}, 0) = -(1+\alpha) \exp\left(-\int (V^{(1)}  + (1-\alpha) z^{(1)} )^{-1} dz^{(1)}\right) \\
&\qquad \times C  \int \exp \left(-\alpha \int (V^{(1)}  + (1-\alpha) z^{(1)} )^{-1} dz^{(1)}\right)  (V^{(1)}  + (1-\alpha) z^{(1)} )^{-2} dz^{(1)}.
\eal
\]From these, we see that $H^2$ can only have one zero at the origin unless it is identically $0$ and the same holds for $W$. Note that $\alpha \approx -2$.
\medskip

  This observation seems to indicate a potential instability in computer calculations at the far field for the SSS.  Since $W(z)$, the profile of vorticity for the expected SSS decay at order around $|z|^{-1/3}$ when $|z|$ is large, a rounding error in computation may cause it to be regarded as $0$ by the computer. According to the observation, this will force $\omega$ and $h$ to be identically $0$ at the base.

Moreover, it shows that $W$ can not become negative in the first quadrant. Indeed, let $z=z(s)$ be a flow line from \eqref{zflowline} and $W(z(s_0))<0$ for some parameter $s_0>0$.
By \eqref{wrect}, we have
\be
\frac{d W(z(s))}{d s} = (\partial_{z^{(1)}} H^2-W)(z(s))>0
\ee for $s$ in a small neighborhood of $s_0$. This shows $W(z(s)) \le W(z(s_0))<0$ for all $s \le s_0$. Since the vector field $(V^{(1)} + (1-\alpha) z^{(1)}, V^{(2)} + (1-\alpha) z^{(2)})$ forms an acute angle with the positive $z^{(1)}$ axis, the flow line can be extended indefinitely backward, we see its backward limit will ether hit the origin or the coordinate axis where $W$ is either $0$ or positive.  This is a contradiction.

However in the file solu.w\{1, 1\} mentioned above, one can find several places where $W$ is negative, including mesh point $(3, 709)$ where $W \approx -2.9778 \times 10^{-18}$.  This shows the approximate steady state there is  not an approximation of the SSS with accuracy within $10^{-18}$ at that mesh point. Note that \eqref{w-5}, which is beyond the stated margin of error,  will also lead to a contradiction.
}

\medskip

{\remark  We would like to discuss the reason for the apparently diverging results. The approximate SSS in \cite{CH1} Sec. 5 is constructed by solving the dynamic equation in a finite domain $D$ for a long time by high performance computing and patching together with a "semi-analytic part" outside this domain using the principle terms of the equation. See the equations 6 lines below (5.1) in \cite{CH1}. This equation is missing the main convection term $v \nabla \omega$ which has the same scaling order as $\theta_x$, which is roughly $r^{-2/3}$. Therefore it is the principal order for the computing error.
As stated, the approximate solution inside $D$ is highly accurate. However, since the convection term is missing for the approximate solution outside of $D$, it is unlikely the two solutions match at the boundary $\partial D$ and the error outside $D$ is largely unknown since $v$ was not computed there.
Even if one finds two exact solutions inside and outside the finite domain with the same boundary value, one may not be able to patch them into a global approximate solution, even less a global exact solution.  One example is the Laplace equation in $\mathbb{R}^3$. Let $D=B(0, 1)$ be the unit ball. Then $f=1$ is a solution in $D$ and $g=1/|x|$ is a solution in $D^c$. It is clear that no global solution of the Laplace equation can be made by patching these two since Liouville theorem shows such global solution must be $0$. It is also clear that by patching them, the high order derivatives near the boundary will not be small, making it less likely to be an approximate solution in $C^k$ sense.}

Finally we consider the case $-2(\beta/\alpha)+(1/\alpha)-1<0$.

\begin{proposition}
\lab{prLHextra}
 Suppose $-2(\beta/\alpha)+(1/\alpha)-1<0$. There does not exist asymptotic self-similar blow-up solutions \eqref{LHsc2} with the following properties. (a). The spatial center of the blow-up sequence is at the side boundary point $p_0=(1, 0, 0)$ in the $(r, \theta, x^{(3)})$ coordinate.  (b). The profile and error are nontrivial  local $C^{2, 2, \gamma}_{x, t}$  functions; (c). $|V^r(y)|$ and $|V^{(3)}(y)|$ are sub-linear near infinity and $|\nabla (V^r, V^{(3)})|$ is bounded ; (d).  $V^{(3)}$ is odd in the vertical variable.
\end{proposition}
 \proof
 If $-2(\beta/\alpha)+(1/\alpha)-1<0$, after takeing $k \to \infty$, we see that the term
$\frac{( \tilde v^\theta_k)^2}{Q_k^{[(1/\alpha)-1]} r}$ in \eqref{eqtvro} vanishes. By the theorem, which also works in the case $\alpha<0$ and $Q_k \to 0$, the limit of $(v^r_k, v^{(3)}_k)$, denoted by $(v^r_k, v^{(3)}_k)$, is a solution of the 2 dimensional Euler equation on the half plane. We can than apply the argument in Corollary \ref{cornossb} (2) to conclude the proposition is true.
 \qed

{\remark
\lab{rmk2}
 From the proof of  Proposition \ref{prLH} part (ii), it is clear that condition
\eqref{LHsc2} can be relaxed to:

(a). for some $C_1, C_2>0$,
\be
\al
&\frac{C_1}{(T_0-t)^{\alpha+\beta}} \Theta \left( \frac{x}{(T_0-t)^{1-\alpha}} \right) \le v^\theta \le \frac{C_2}{(T_0-t)^{\alpha+\beta}} \Theta \left( \frac{x}{(T_0-t)^{1-\alpha}} \right)  \\
&\frac{C_1}{(T_0-t)^\alpha}  V^r\left( \frac{x}{(T_0-t)^{1-\alpha}} \right) \le v^r \le  \frac{C_2}{(T_0-t)^\alpha}  V^r\left( \frac{x}{(T_0-t)^{1-\alpha}} \right)  \\
&\frac{C_1}{(T_0-t)^\alpha}  V^{(3)} \left( \frac{x}{(T_0-t)^{1-\alpha}} \right) \le v^{(3)} \le  \frac{C_2}{(T_0-t)^\alpha}  V^{(3)} \left( \frac{x}{(T_0-t)^{1-\alpha}} \right);
\eal
\ee

(b). The scaled functions $Q^\beta_k \tilde v^\theta_k$, $\tilde v^r_k$ and $\tilde v^{(3)}_k$ converge in local $C^{2, 1, \gamma^-}_{x, t}$ sense.
}

\medskip

We finish by mentioning a generalization of Corollary \ref{cornossb} to a class of solutions which are called quasi self-similar.

\begin{defn} A solution $v$ to \eqref{EL} is called  {\it quasi self-similar} (QSS) with quasi profile $V$  if
\[
v=v(x, t) = \frac{1}{(T_0-t)^\alpha} V \left( \frac{x}{(T_0-t)^{1-\alpha}}, t \right)
\] for some bounded vector field $V$ and some $\alpha \in (0, 1)$.
\end{defn}

 Assuming the quasi profile $V(\cdot, \cdot) \in C^{2, 1, \gamma}_{x, t}$, then the conclusions of Corollary \ref{cornossb} still hold. The proof is almost verbatim since the scaled solution
 \[
 \tilde v_k=\frac{1}{(\, 1 - \tilde t \, )^{\alpha} } V\left(\frac{\tilde x + x_k Q_k^{(1-\alpha)/\alpha}}{(1-\tilde t)^{1-\alpha}}, \,  Q_k^{-1/\alpha} \tilde t + t_k \right)
 \]are uniformly  bounded in $C^{2, \gamma}_{\tilde x}$ norm for all $k$ and $\tilde t \le 0$.

\section*{Acknowledgment} The author is grateful to the support of the Simons
Foundation through grant No. 710364.
He would also like to thank Professors Tristan Buckmaster, Jiajie Chen, Tom Hou for useful communications and Professors Hongjie Dong, Zijin Li,  Xinghong Pan, Na Zhao and Drs Xin Yang and Chulan Zeng for discussions. He is especially grateful to Prof. Jiajie Chen for pointing out an error in Proposition 2.1 (ii), i.e. $2 \beta=1-\alpha$ case,  in version 3 of the paper in the arxiv and many discussions \cite{Chj}.

\end{document}